\def\pidpreprintmode{1}
\documentclass[10pt]{article}
\newif\ifpidpreprint\pidpreprintfalse
\ifdefined\pidpreprintmode\pidpreprinttrue\fi
\ifpidpreprint\usepackage[preprint]{tmlr}\else\usepackage{tmlr}\fi
\ifpidpreprint\fi
\fancypagestyle{pidfirst}{\fancyhead{}\renewcommand{\headrulewidth}{0pt}}
\usepackage[utf8]{inputenc}
\usepackage{microtype,amsmath,amssymb,amsthm,booktabs,tabularx,array}
\usepackage{hyperref}
\usepackage{orcidlink}
\hypersetup{colorlinks=true,linkcolor=blue!45!black,citecolor=blue!45!black,urlcolor=blue!45!black,
  pdftitle={Your PID Is Already a Neural Network: Learning Its Gains with Only Three Weights and No Teacher}}
\ifpidpreprint\hypersetup{pdfauthor={Felipe de Jes\'us F\'elix-Arredondo, Nezih Nieto-Guti\'errez, Carlos Astengo-Noguez and Alejandro Ucan-Puc}}\fi
\usepackage{url}
\usepackage{xcolor}
\definecolor{pidorange}{HTML}{EB6834}\definecolor{recblue}{HTML}{2A78D6}\definecolor{refgray}{HTML}{52514E}
\usepackage{colortbl}
\definecolor{passbg}{HTML}{DDEFDD}\definecolor{failbg}{HTML}{F8DCD7}\definecolor{notebg}{HTML}{FBEFD3}
\newcommand{\cpass}{\cellcolor{passbg}}\newcommand{\cfail}{\cellcolor{failbg}}

\definecolor{passmark}{HTML}{5BA867}\definecolor{notemark}{HTML}{E0A033}
\newcommand{\mpass}{\textcolor{passmark}{\rule[-0.15ex]{0.8em}{0.8em}}\hspace{0.45em}}
\newcommand{\mnote}{\textcolor{notemark}{\rule[-0.15ex]{0.8em}{0.8em}}\hspace{0.45em}}
\usepackage{etoolbox}
\AtBeginEnvironment{table}{\setlength{\aboverulesep}{0pt}\setlength{\belowrulesep}{0pt}%
  \setlength{\abovecaptionskip}{0pt}\setlength{\belowcaptionskip}{5pt}}
\newcommand{\rowgap}{\addlinespace[1.5pt]}
\usepackage{enumitem}
\usepackage{graphicx}
\graphicspath{{figures/}}
\usepackage{tikz}
\usetikzlibrary{positioning,arrows.meta,calc}
\usepackage[section]{placeins}
\usepackage{needspace}
\setlist{nosep}
\usepackage[capitalise,noabbrev]{cleveref}

\edef\pidjob{\jobname}\edef\pidfinaljob{\detokenize{PIDAsANet-final}}
\ifx\pidjob\pidfinaljob\def\pidfinal{1}\else\def\pidfinal{0}\fi
\newcommand{\pending}[1]{\ifnum\pidfinal=1
  \PackageError{PIDAsANet}{Unresolved \string\pending: #1}{Every pending result must be
  replaced by a generated number before a final build.}%
  \else\fcolorbox{red}{red!6}{\textcolor{red}{\sffamily\small pending: #1}}\fi}
\newtheorem{theorem}{Theorem}
\newtheorem{proposition}[theorem]{Proposition}
\newtheorem{corollary}[theorem]{Corollary}
\theoremstyle{remark}\newtheorem{remark}[theorem]{Remark}
\newcommand{\E}{\mathbb E}
\newcommand{\R}{\mathbb R}
\newcommand{\PID}{\mathrm{PID}}
\newcommand{\T}{\mathsf T}
\newcommand{\Id}{\mathrm{Id}}
\newcommand{\EoneCases}{10{,}000}
\newcommand{\EoneMax}{6.2\times10^{-13}}

\newcommand{\EtwoReps}{1000}
\newcommand{\EtwoMaxAbsZ}{1.93}

\newcommand{\EtwoReadoutRate}{0.928203}

\newcommand{\EtwoNEMaxCorr}{6.4\times10^{-9}}
\newcommand{\LpidStar}{3671/2750}
\newcommand{\GapRelPid}{5.88\%}
\newcommand{\GapRelPidBias}{4.00\%}
\newcommand{\EthreeRealizableMax}{3.4\times10^{-10}}
\newcommand{\EthreeReps}{1000}
\newcommand{\EthreebDivPid}{0.52\%}
\newcommand{\EthreebDivPidBias}{0.30\%}

\newcommand{\EfourEffortShareLinear}{4\%}
\newcommand{\EfourIaeEnvLinear}{0.271}
\newcommand{\EfourIaeWarmLinear}{0.275}
\newcommand{\EfourWarmLinear}{$-$3.51\%}
\newcommand{\EfourWarmUpLinear}{$-$2.49\%}
\newcommand{\EfourBoundedLinear}{$-$3.70\%}
\newcommand{\EfourBoundedLoLinear}{$-$4.82\%}
\newcommand{\EfourEnvEqLinear}{$-$14.27\%}
\newcommand{\EfourWarmShiftLinear}{$-$2.01\%}
\newcommand{\EfourRhoRatioLinear}{1.29}
\newcommand{\EfourFirstOrderLinear}{0.997}
\newcommand{\EfourEffortShareDuffing}{3\%}
\newcommand{\EfourIaeEnvDuffing}{0.272}
\newcommand{\EfourIaeWarmDuffing}{0.273}
\newcommand{\EfourWarmDuffing}{$-$3.45\%}
\newcommand{\EfourWarmUpDuffing}{$-$2.54\%}
\newcommand{\EfourBoundedDuffing}{$+$14.03\%}

\newcommand{\EfourEnvEqDuffing}{$-$7.36\%}
\newcommand{\EfourWarmShiftDuffing}{$-$2.30\%}
\newcommand{\EfourRhoRatioDuffing}{1.00}
\newcommand{\EfourFirstOrderDuffing}{0.999}
\newcommand{\EfourEffortShareHammerstein}{1\%}
\newcommand{\EfourIaeEnvHammerstein}{0.057}
\newcommand{\EfourIaeWarmHammerstein}{0.074}
\newcommand{\EfourWarmHammerstein}{$-$1.90\%}
\newcommand{\EfourWarmUpHammerstein}{$-$1.38\%}
\newcommand{\EfourBoundedHammerstein}{$+$3.00\%}

\newcommand{\EfourEnvEqHammerstein}{$-$2.72\%}
\newcommand{\EfourWarmShiftHammerstein}{$-$2.56\%}
\newcommand{\EfourRhoRatioHammerstein}{0.34}
\newcommand{\EfourFirstOrderHammerstein}{0.992}
\newcommand{\EfourBootstrap}{10{,}000}
\newcommand{\EfourAttrPidLinear}{$-$1.07\%}
\newcommand{\EfourAttrResPidLinear}{$-$2.46\%}
\newcommand{\EfourAttrResPidUpLinear}{$-$1.60\%}

\newcommand{\EfourAttrShareLinear}{30\%}
\newcommand{\EfourAttrPidDuffing}{$-$0.75\%}
\newcommand{\EfourAttrResPidDuffing}{$-$2.73\%}
\newcommand{\EfourAttrResPidUpDuffing}{$-$1.95\%}

\newcommand{\EfourAttrShareDuffing}{22\%}

\newcommand{\EfourAttrResPidHammerstein}{$-$0.25\%}

\newcommand{\EfourAttrShareHammerstein}{87\%}

\newcommand{\EfiveLoops}{504}
\newcommand{\EfiveSets}{21}
\newcommand{\EfiveMaxTwoState}{2\times10^{-13}}
\newcommand{\EfiveMaxAll}{4\times10^{-12}}

\newcommand{\EfiveLearnErr}{4\times10^{-10}}
\newcommand{\EfiveLearnLoop}{3\times10^{-10}}
\newcommand{\EfiveItersGain}{156}
\newcommand{\EfiveItersNet}{35{,}668}
\newcommand{\EfiveCond}{1{,}600}

\newcommand{\LyapMassPid}{88.4\%}
\newcommand{\LyapMassStar}{99.99\%}
\newcommand{\LyapDivPid}{-3}
\newcommand{\LyapDivStar}{-7}
\newcommand{\LyapEfiveRel}{5\times10^{-16}}

\title{Your PID Is Already a Neural Network:\\ Learning Its Gains with Only Three Weights\\ and No Teacher}
\author{\name Felipe de Jes\'us F\'elix-Arredondo~\orcidlink{0009-0006-6894-5694} \email ffelix3ro@gmail.com\\
      \addr Healthen Group, Mexico
      \AND
      \name Nezih Nieto-Guti\'errez~\orcidlink{0000-0001-8512-005X} \email nezihniegu@tec.mx\\
      \addr School of Engineering and Sciences, Tecnol\'ogico de Monterrey, Mexico
      \AND
      \name Carlos Astengo-Noguez~\orcidlink{0009-0001-2433-3817} \email castengo@tec.mx\\
      \addr School of Engineering and Sciences, Tecnol\'ogico de Monterrey, Mexico
      \AND
      \name Alejandro Ucan-Puc~\orcidlink{0000-0002-0037-9394} \email alejandro.ucan-puc@tec.mx\\
      \addr School of Engineering and Sciences, Tecnol\'ogico de Monterrey, Mexico}
\def\month{MM}\def\year{YYYY}\def\openreview{\url{https://openreview.net/forum?id=XXXX}}
\begin{document}
\maketitle\ifpidpreprint\thispagestyle{pidfirst}\fi
\begin{abstract}
A discrete PID controller is exactly a recurrent neural network, a classical state-space fact: a fixed linear memory of the running integral and the last error, read out by three trainable weights that map one-to-one to the PID gains. Networks with bounded hidden states, such as LSTMs of any width, only approximate it. Trained only on a one-step physical cost (the expected cost of one transition from a fixed state distribution), with no demonstrations, teacher controller, or known gains, the network solves a strongly convex least-squares problem, so gradient descent converges from any initialization to the unique one-step optimal PID gains; the limit network is that PID, with the same action on every error sequence. Partitioned regression then shows what extra capacity buys. Orthogonal residual features keep these gains and improve on the PID if and only if they correlate with the part of the optimal action outside the PID span. On a nonlinear second-order example the gap is $5.88\%$ of the optimal PID cost, and $4.00\%$ against PID with a bias; no controller can improve on PID by more on this objective. A residual started at a tuned PID begins exactly at that PID; at the one-step optimal PID, a fixed-feature residual's first gradient step lowers the cost if and only if its features can improve on the PID. Experiments confirm the realization, recovery, and gap results to numerical precision. On held-out plants of three families, under a protocol frozen before the test data existed, a residual warm-started at a heavily tuned PID and trained with it lowers the long-horizon cost by 3.5\%, 3.5\%, and 1.9\%, with family-wise 95\% bounds excluding zero, also under distribution shift; a post hoc comparison with the PID tuned alone suggests that most of it comes from the recurrent part in two families. The primary endpoint, a gain of at least 5\% from a generic PID start, is not met in any family.
\end{abstract}

\section{Introduction}
Most feedback loops in practice are still closed by proportional--integral--derivative (PID) controllers~\citep{astrom2001}, and learned controllers are often proposed to replace them. A discrete PID compresses the error history into three features: the present error, the accumulated error, and a finite difference. That compression is already a recurrent computation. A discrete PID is a recurrent neural network with a fixed linear memory, which holds the running integral and the last error, and a readout whose three weights map one-to-one to the gains (\cref{fig:pidrnn}), exactly, for every gain, error sequence, and sampling period, and also in continuous time once the derivative is filtered. The representation is classical: linear recurrent networks are the linear state-space systems of control theory, and this one is a standard state-space form of PID~\citep{astrom2006}, so we claim no novelty for it. Nor is it a linearized LSTM. A standard LSTM has a bounded output, so it cannot reproduce integral action over all times (\cref{prop:lstm}), although it can learn the integral over the horizons it is trained on. A learned controller can therefore include PID exactly instead of approximating it, and the question becomes what such a network learns when PID is part of it, and when capacity beyond PID pays off.

Our main result is that the network finds the optimal PID gains by itself. It is trained only on the cost of plant transitions, with no demonstrations, no teacher controller, and no knowledge of which gains are good, and gradient descent converges from every initialization to the unique gains that are optimal for that cost; the converged network is then that PID. The cost is deliberately simple: the expected cost of one plant transition, with penalties on the next state and on the action, for states drawn from a fixed distribution. It is a quadratic criterion of the kind that underlies classical self-tuning control, where a PID structure can be derived from it~\citep{clarke1975,cameron1983}. We use it because in this setting learning the gains and the value of added capacity separate exactly, so that every prediction can be checked to numerical precision, and because it also describes, for a residual with fixed features, the first gradient step from the one-step optimal PID, the counterpart of the warm start used over long horizons (\cref{sec:warm}).

Because the memory is fixed and the gains sit in the readout, this objective is quadratic in the readout parameters even when the drift of the plant is nonlinear, as long as the action enters affinely, which makes the convergence global and geometric, with an exact finite-sample error law (\cref{sec:realization,sec:learning}). To measure what extra capacity buys, we add fixed nonlinear features to the readout and orthogonalize them against the PID features. The PID block then keeps its optimal gains, and the added block lowers the cost by a closed-form amount, which is positive if and only if the features correlate with the part of the optimal action that no PID can represent (\cref{sec:improvement}). On a nonlinear second-order plant this gap is 5.88\% of the optimal PID cost, a size set by the strength of the nonlinearity chosen there. On classical PID benchmarks, from a mass-spring system in its three damping regimes to an inverted pendulum, coupled tanks, and a thermal chamber with dead time, the network and the textbook PID produce the same closed loop to rounding, and the network trained on the physical cost learns the one-step optimal gains (\cref{sec:results}).

None of these ingredients is new on its own: besides the realization, gain recovery is least squares, and coefficient preservation is the Frisch--Waugh--Lovell theorem of partitioned regression~\citep{frisch1933,lovell1963}. What we add is a setting in which each prediction is exact and all but (iv) are checked numerically: (i)~readout coordinates for the recurrent realization and their conditioning, which explain why a step size that is safe for the gains can diverge for the recurrent weights; (ii)~teacher-free identification of the one-step optimal gains, with an exact finite-sample error law; (iii)~a closed-form gap on a nonlinear plant; (iv)~the observation that, at the one-step optimal PID, the first gradient step of a zero-initialized residual with fixed features lowers the cost if and only if that gap is positive; and (v)~a long-horizon evaluation under a protocol frozen before its test data existed.

Over long horizons the visited states depend on the controller, and the theory says only that a residual initialized at zero starts exactly at its PID and what its first gradient step does (\cref{sec:long}). The long-horizon experiment E4 tests the rest on three benchmark families, against the PID envelope, the best PID found by a search that used 650 to 1,300 times the budget of one learned training run. A residual whose PID part starts from generic gains does not improve on the envelope by the 5\% of the primary endpoint in any family, and is worse than the envelope on the two nonlinear families. One whose PID part starts at the envelope improves on it in all three, by 1.9--3.5\%, on held-out plants and under distribution shift, and inherits the cost of that tuning. A post hoc comparison with the PID tuned alone from the same start suggests that most of this gain comes from the recurrent part in the two oscillator families and from further tuning of the PID in the Hammerstein family (\cref{sec:attribution}). We do not claim that neural PID, residual control, or differentiation through feedback is new, nor that a one-step optimum is optimal or stabilizing over an infinite horizon.

\Cref{sec:related} reviews related work. \Cref{sec:realization,sec:learning,sec:improvement} develop the one-step theory, \cref{sec:long} moves to long horizons and the warm start, \cref{sec:experiments,sec:results} describe and report the experiments, and \cref{sec:conclusions} discusses the limitations. The appendices give proofs, protocols, full results, and notation.

\section{Related work and scope of the contribution}\label{sec:related}
Recurrent neural networks have already been used to tune PID gains online. \citet{gunther} evaluate general dynamic neural networks for that purpose, and PID neural networks~\citep{shu2000} build hidden neurons that compute proportional, integral and derivative signals. Our fixed-memory realization instead makes the PID gains readout parameters of the network. This is an elementary observation; its use here is to make PID a submodel of a larger network and to fix the coordinates in which learning is analyzed.

Tuning PID gains from data has a long history in control. Self-tuning controllers choose the control law that minimizes a quadratic cost of the future output, the set point, and the control effort, with model parameters estimated recursively from data~\citep{clarke1975}, and \citet{cameron1983} derive a self-tuning controller with a PID structure in this way. Our one-step objective is a cost of the same kind. The difference is that we fix a distribution of states and minimize the expected cost directly over the gains, which makes recovery a regression with an exact finite-sample law. Virtual reference feedback tuning~\citep{campi2002} identifies controller parameters by least squares from a single data set, without a plant model but with a reference model that specifies the desired closed loop, and iterative feedback tuning~\citep{hjalmarsson1998} estimates the gradient of a closed-loop cost from experiments on the plant. Our gradient result needs the input gain of the plant (assumption (A5) in \cref{sec:learning}), or else costs of alternative actions from the same reset states and noise (\cref{app:theory}), so it is closer to model-based design than these data-driven methods.

DiffLoop~\citep{diffloop} differentiates through a feedback loop to optimize PID gains and anti-windup behavior, and DiffTune~\citep{cheng2024difftune} autotunes controller parameters by differentiating through closed-loop simulations more generally, so learning from a physical cost without gain labels is not new either. \citet{sharma} establish policy-gradient convergence guarantees for PID policies in their setting. Our convergence theorem is narrower. It concerns one transition from a fixed distribution of states, where the loss is quadratic in the gains, and it does not replace an analysis of the closed-loop landscape. That landscape, which changes with the policy, is what policy optimization for control studies, from global convergence for the linear quadratic regulator~\citep{fazel2018} to the broader theory surveyed by \citet{hu2023}.

Residual policy learning~\citep{silver2018,johannink2019} has also been combined with PID, including the work of \citet{hynes}. Zero-initialized residuals on top of a base controller are standard in that line of work, so we claim no priority for the warm start itself. What we add is an exact account of what a residual can contribute in the one-step setting: an orthogonal decomposition gives both coefficient preservation and a closed-form improvement. These identities are classical; we use them to analyze the first training step at the one-step optimal PID (\cref{cor:warm}), the counterpart of the long-horizon warm start whose first step (\cref{prop:warm}) \cref{sec:results} tests.

Fractional-order PID controllers enlarge the structure differently, with integral and derivative operators of arbitrary order; \citet{zarate2023} report improved transients with such controllers on an inverted pendulum. They remain fixed linear operators of the error history, whereas the residual studied here adds learned nonlinear features while keeping the PID gains as an explicit submodel.

Conventional LSTM cells~\citep{lstm} suggest reading the cell state as a gated memory, but the identity memory path and hard gate values that an exact integral needs are modifications of the standard cell. LSTMs can learn to count in practice~\citep{gers2001,weiss2018}, but a bounded output cannot track an unbounded integral globally (\cref{prop:lstm}). Linear recurrent memories such as the linear recurrent unit (LRU)~\citep{orvieto2023} retain the near-identity dynamics that the exact realization uses.

\section{Discrete PID as a recurrent network}\label{sec:realization}
The realizations in this section are the classical state-space and velocity (incremental) forms of the discrete PID law~\citep{astrom2006}. We restate them for two reasons. The realization makes PID a submodel of a larger recurrent network, which the residual extension and the warm start of \cref{sec:long} rely on. It also fixes the coordinates of the network's weights, in which, as \cref{sec:learning} shows, training the gains behaves differently than in the gains themselves.

Let $\Delta>0$ be the sampling period and $e_t=r_t-y_t$ the error between reference and plant output at step $t\geq0$; $I_{-1}$ and $e_{-1}$ are fixed initial values. The discrete PID approximates the integral of the error by a running sum and its derivative by a backward difference:
\begin{equation}\label{eq:pid}
 I_t=I_{t-1}+\Delta e_t,\qquad D_t=\frac{e_t-e_{t-1}}{\Delta},\qquad
 u_t=k^\T\phi_t,\quad \phi_t=(e_t,I_t,D_t)^\T,
\end{equation}
where $k=(K_p,K_i,K_d)^\T\in\R^3$ are the gains and $\phi_t$ the PID features. The action is linear in the gains, which the learning results below rely on. The analytic theorems use this unrestricted gain class; gain constraints, filtering, and saturation enter in \cref{sec:long}. \Cref{app:notation} lists the notation.

\begin{figure}[!b]
\centering
\begin{tikzpicture}[>={Stealth[length=2mm]}, every node/.style={font=\small},
  blk/.style={draw, rounded corners=2pt, align=center, inner sep=4pt, minimum height=9mm},
  mem/.style={blk, draw=refgray, fill=refgray!8},
  pidb/.style={blk, draw=pidorange, fill=pidorange!10, line width=0.6pt},
  resb/.style={blk, draw=recblue, fill=recblue!8, line width=0.6pt},
  sm/.style={draw, circle, inner sep=1pt}]
\node (in) {$e_t$};
\node[mem, right=16mm of in] (cell) {fixed linear memory\\[1pt] $s_t=A_ss_{t-1}+b_se_t$};
\node[mem, right=12mm of cell] (dly) {delay};
\node[pidb, below=9mm of dly] (read) {learnable readout $(c,c_0)$\\[1pt] $c^\T s_{t-1}+c_0e_t$};
\node[sm, right=9mm of read] (plus) {$+$};
\node[right=7mm of plus] (out) {$u_t$};
\node[resb, dashed, below=5mm of read] (res) {optional residual $\beta^\T\psi$ (\cref{th:gap})};
\coordinate (tap) at ($(in.east)+(6mm,0)$);
\draw[->] (in) -- (cell);
\draw[->] (cell) -- node[above] {$s_t$} (dly);
\draw[->] (dly.north) -- ++(0,5mm) -| node[pos=0.25, above] {$s_{t-1}$} (cell.north);
\draw[->] (dly.south) -- node[right] {$s_{t-1}$} (read.north);
\draw[->] (tap) |- (read.west);
\draw[->, dashed, recblue] (tap) |- (res.west);
\draw[->] (read) -- (plus);
\draw[->, dashed, recblue] (res.east) -| (plus.south);
\draw[->] (plus) -- (out);
\fill (tap) circle (1pt);
\end{tikzpicture}
\caption{PID as a recurrent submodel (\cref{th:real}). The memory is fixed and linear; all PID gains live in the readout $(c,c_0)$, which maps bijectively to $(K_p,K_i,K_d)$. The dashed branch is the residual extension of \cref{th:gap}; with $\beta=0$ the controller is exactly the PID.}
\label{fig:pidrnn}
\end{figure}
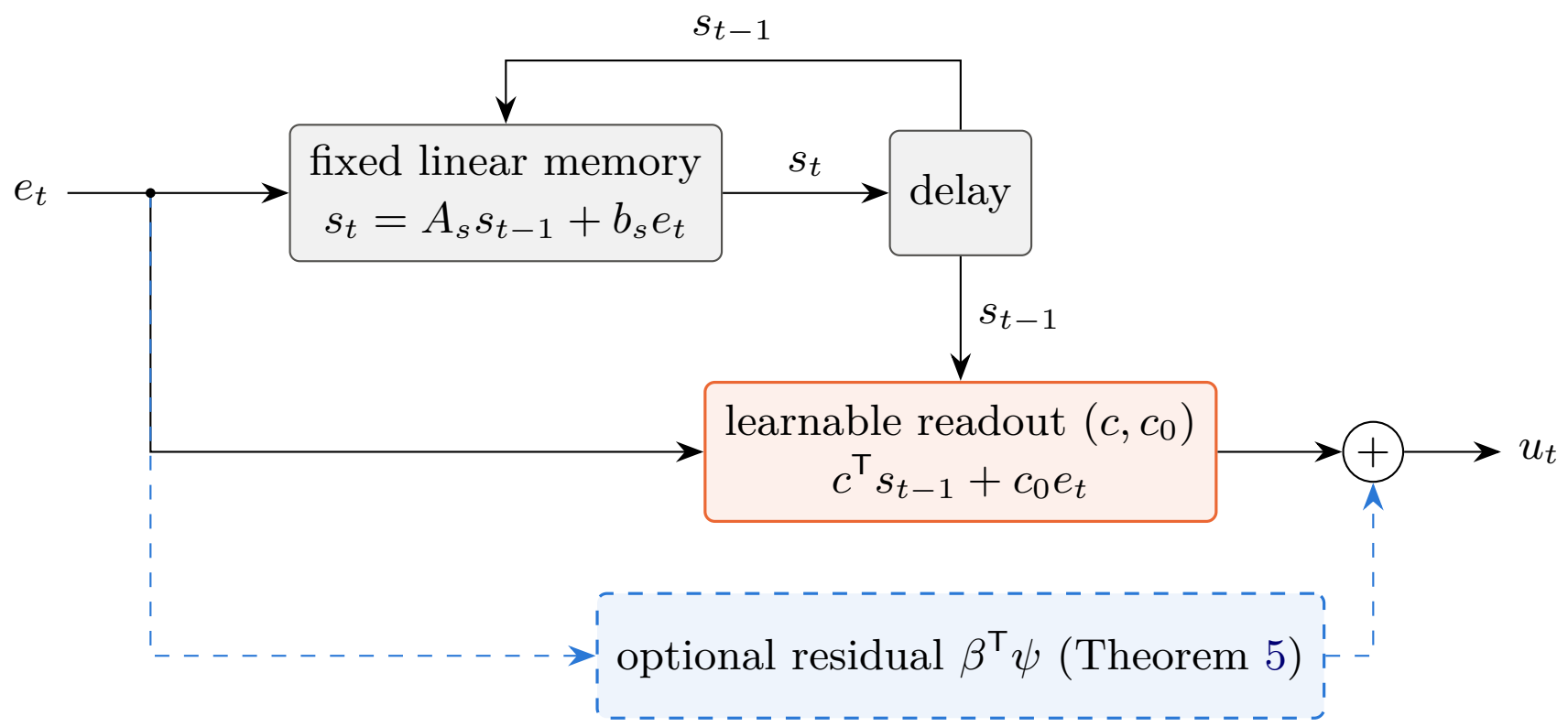

\begin{theorem}[Exact two-state realization]\label{th:real}
Every controller in \eqref{eq:pid} is exactly represented, for every real error sequence and every finite time, by
\begin{equation}
 s_t=A_s s_{t-1}+b_s e_t,\qquad
 A_s=\begin{bmatrix}1&0\\0&0\end{bmatrix},\quad b_s=\begin{bmatrix}\Delta\\1\end{bmatrix},
 \quad s_{-1}=\begin{bmatrix}I_{-1}\\e_{-1}\end{bmatrix},
\end{equation}
with readout
\begin{equation}\label{eq:readout}
 u_t=c^\T s_{t-1}+c_0e_t,\quad
 c=\begin{bmatrix}K_i\\-K_d/\Delta\end{bmatrix},\quad
 c_0=K_p+\Delta K_i+K_d/\Delta.
\end{equation}
The readout parameterization is bijective, with inverse $K_i=c_1$, $K_d=-\Delta c_2$, and $K_p=c_0-\Delta c_1+c_2$.
\end{theorem}
\begin{proof}
Induction gives $s_t=(I_t,e_t)^\T$: the first coordinate accumulates $\Delta e_t$, and the second is overwritten by $e_t$. Expanding the readout and regrouping,
\begin{equation*}
 c^\T s_{t-1}+c_0e_t=K_iI_{t-1}-\frac{K_d}{\Delta}e_{t-1}+\Big(K_p+\Delta K_i+\frac{K_d}{\Delta}\Big)e_t
 =K_pe_t+K_i(I_{t-1}+\Delta e_t)+K_d\frac{e_t-e_{t-1}}{\Delta},
\end{equation*}
which is $K_pe_t+K_iI_t+K_dD_t$, the action \eqref{eq:pid}. The map from $k$ to $(c,c_0)$ is linear, and solving its three equations for $k$ gives the inverse.
\end{proof}
\Cref{fig:pidrnn} shows the realization. We use this two-state form as the main realization because its memory is minimal, two fixed linear functions of the error history, and all three gains sit in the readout. The price is the feedthrough term $c_0e_t$, which a standard recurrent readout $u_t=Wh_t$ lacks. \Cref{app:theory} gives a current-state (three-state) realization without feedthrough, the incremental form, and a gated linear-memory variant; all realize the same operator; E1 checks them on error sequences and E5 in closed loop on classical plants (\cref{sec:experiments}).

\begin{proposition}[Bounded-output recurrent readouts]\label{prop:lstm}
Let $u_t=Wh_t+b_o$ with $W\in\R^{1\times n_h}$ and $b_o\in\R$, where $|h_{t,j}|\leq1$ for all $t$ and $j$. This holds for a conventional LSTM of any width and with any weights, whose hidden state $h_t=o_t\odot\tanh(c^{\rm cell}_t)$ is the elementwise product of an output gate in $(0,1)$ and the $\tanh$ of the cell state. Then $|u_t|\leq\sum_j|W_j|+|b_o|$ for all $t$. Consequently, no such network realizes a PID \eqref{eq:pid} with $K_i\neq0$ over all times, already for the bounded error $e_t\equiv1$. For $K_i=0$ and $(K_p,K_d)\neq(0,0)$, exact realization fails on unbounded error sequences.
\end{proposition}
\begin{proof}
Since $|h_{t,j}|\leq1$, $|u_t|\leq\sum_j|W_j||h_{t,j}|+|b_o|\leq\sum_j|W_j|+|b_o|$. For $e_t\equiv1$ and $t\geq1$, \eqref{eq:pid} gives $u_t=K_p+K_i(I_{-1}+\Delta(t+1))$, which is unbounded when $K_i\neq0$. For $K_i=0$ and $e_t=t^2$, $u_t=K_pt^2+K_d(2t-1)/\Delta$ for $t\geq1$, which is unbounded unless $K_p=K_d=0$.
\end{proof}
The realization therefore keeps the integral in an explicit linear memory. A bounded network can still be added on top of it, as the residual of \cref{sec:warm} is, because the PID part then supplies the integral action. The bound concerns exact, unsaturated realization through a bounded hidden state; \cref{app:theory} discusses PD controllers, cell-state readouts, actuator saturation, and why PID and its recurrent realization have the same expected trajectory cost.

\paragraph{Sampling period and continuous time.} Nothing above depends on the sampling period: \cref{th:real} holds for every $\Delta>0$, E1 checks it for periods from $10^{-3}$ to $1$, and E5 at the periods of classical benchmarks. The period enters only through the readout, whose derivative weight scales as $1/\Delta$, which is why the readout coordinates become ill-conditioned as $\Delta\to0$ (\cref{cor:readout}). Continuous-time PID, with the filtered derivative it needs to be realizable, is likewise a continuous-time linear recurrent network, with a fixed memory (the integral and the filtered error) and a readout that maps one-to-one to the gains; sampled with the error held over each period, it is a discrete fixed linear memory with the same readout (\cref{prop:continuous}). The learning results of \cref{sec:learning,sec:improvement} need only an action that is linear in the gains and computed from a fixed memory, so they apply to that sampled controller as well. They concern sampled transitions, as any digital implementation does; learning in continuous time is outside our scope.

\section{Learning optimal gains without a PID teacher}\label{sec:learning}
\Cref{sec:realization} shows that the network can take the form of any PID, but not which gains training reaches; fitting the readout to the actions of a known PID would presuppose them. This section shows that the gains need not be known in advance. For a declared cost and distribution of states, and under full-rank excitation, the one-step optimal PID gains are the unique solution of an optimization problem built only from plant transitions and their costs (\cref{th:recovery}). Gradient descent finds them without a teacher controller, and because the readout is a bijection (\cref{th:real}), the converged network is that PID.

\subsection{Physical objective and sampling assumptions}
Let $z\in\R^n$ be the augmented state: the observed plant state together with the controller memory, such as the integral $I$, so that the next state depends only on $z$, the action, and noise. For the PID, $z$ contains the current error $e_t$, the integral $I_t$ already updated with it, and the previous error $e_{t-1}$ or, equivalently, the difference $e_t-e_{t-1}$ (as in \cref{sec:example}), so that the features $\phi_t$ of \eqref{eq:pid} are a function of $z$. States are drawn from a design distribution $\nu$, and the PID features are functions $\phi(z)$ of $z$. The next state $z^+$ follows a transition in which the scalar action $u$ enters affinely,
\begin{equation}\label{eq:plant}
 z^+=a(z)+bu+\xi,
\end{equation}
where the drift $a(z)$ is the next state under zero action, $b\in\R^n$ is a constant input vector, and $\xi$ is transition noise. For a linear plant $a(z)=Az$ for a matrix $A$; \cref{sec:example} constructs $a$ and $b$ for a nonlinear plant. We assume that the desired next state is zero, possibly after an explicit change to reference coordinates. For $H\succeq0$ and $R>0$, selected before optimization, define
\begin{equation}\label{eq:physical}
 \mathcal L(u)=\E_{z\sim\nu,\xi}\big[(z^+)^\T H z^++Ru^2\big],
 \qquad M=R+b^{\T}Hb>0.
\end{equation}
This criterion charges the next state and the action, but not the states that follow, as in the one-step-ahead costs of self-tuning control (\cref{sec:related}). The results of this section rely on the following assumptions. Each fails in some practical setting, and E4 includes cases that violate (A1) and (A3).
\begin{enumerate}[label=(A\arabic*)]
\item \emph{Control-affine transition.} The action enters \eqref{eq:plant} affinely through a constant $b\in\R^n$. A state-dependent $b(z)$ turns the regressions below into weighted least squares, with weight $M(z)=R+b(z)^\T Hb(z)$ (\cref{app:theory}); a nonlinear input map, such as the Hammerstein family of \cref{sec:experiments}, destroys the quadratic structure.
\item \emph{Noise.} $\E[\xi\mid z,u]=0$, and $\E[\xi^\T H\xi\mid z,u]$ is finite and independent of $u$.
\item \emph{Fixed design distribution.} The distribution $\nu$ is fixed independently of the learned policy: a reset-state or offline design distribution, not the distribution that the current controller induces.
\item \emph{Deterministic policies.} The policy is a measurable function $u=\pi(z)$. By \eqref{eq:square} below, a randomized policy has the loss of its conditional-mean policy plus $M\,\E[\operatorname{Var}(u\mid z)]\geq0$, so this restriction loses no optimality.
\item \emph{Gradient access.} Training evaluates \eqref{eq:physical} on observed transitions and never observes optimal actions or PID coefficients, so it is teacher-free. It is not model-free: an exact gradient needs $\partial z^+/\partial u=b$, that is, a differentiable transition model; \cref{app:theory} gives an exact derivative-free alternative.
\end{enumerate}
We call \eqref{eq:physical} a physical objective because it is computed from transitions and their costs alone, with no teacher actions. When $z$ contains controller memory, such as the integral $I$, the weight $H$ can penalize that memory too (\cref{sec:example}).

Under (A1)--(A4) the loss of a policy $u=u(z)$ takes a form that shows its minimizer. Write $a=a(z)$ and substitute \eqref{eq:plant} into \eqref{eq:physical}. The terms linear in $\xi$ vanish by (A2), and
\begin{equation*}
 \mathcal L(u)=\E\big[a^\T Ha+2u\,b^\T Ha+Mu^2\big]+\E[\xi^\T H\xi].
\end{equation*}
For each $z$ the bracket is a quadratic in the scalar $u$ with leading coefficient $M>0$. Completing the square, $Mu^2+2u\,b^\T Ha=M(u-u_\star)^2-Mu_\star^2$ with $u_\star=-M^{-1}b^\T Ha$, gives
\begin{equation}\label{eq:square}
 u_\star(z)=-M^{-1}b^{\T}Ha(z),\qquad
 \mathcal L(u)=\mathcal L_{\rm irr}+M\E[(u(z)-u_\star(z))^2],
\end{equation}
where the irreducible cost $\mathcal L_{\rm irr}=\E[a^\T Ha-Mu_\star^2+\xi^\T H\xi]$ does not depend on the policy. The one-step optimal action $u_\star(z)$ is therefore the best action at every state, and the loss of any policy is $\mathcal L_{\rm irr}$ plus $M$ times its mean squared distance from $u_\star$. The function $u_\star$ appears only in the analysis; training uses \eqref{eq:physical}.

\subsection{Global coefficient convergence}
Let $\phi(z)\in\R^3$ denote the PID features, assume all moments below are finite, and set
\begin{equation*}
 \Sigma=\E[\phi\phi^\T],\qquad m_\phi=\E[\phi\,u_\star].
\end{equation*}
The condition $\Sigma\succ0$ is the full-rank excitation used below: no nonzero combination $k^\T\phi$ vanishes almost surely under $\nu$, so the three features are linearly independent as functions of the state. \emph{Full-batch population gradient descent} is the iteration $k_{\ell+1}=k_\ell-\eta\nabla\mathcal L(k_\ell)$ on the loss \eqref{eq:physical} itself, with exact expectations instead of sample averages, step size $\eta>0$, and iteration index $\ell$.
\begin{theorem}[Teacher-free optimal PID recovery]\label{th:recovery}
If $\Sigma\succ0$, the physical loss restricted to $u=k^\T\phi$ has the unique minimizer
\begin{equation}\label{eq:kstar}
 k_\star=-M^{-1}\Sigma^{-1}\E[\phi\,b^{\T}Ha(z)]
 =\Sigma^{-1}m_\phi.
\end{equation}
Its Hessian is $2M\Sigma$. Full-batch population gradient descent with
$0<\eta<1/(M\lambda_{\max}(\Sigma))$ satisfies
\begin{equation}\label{eq:rate}
 \|k_\ell-k_\star\|_2\leq q^\ell\|k_0-k_\star\|_2,\quad
 q=\max_j|1-2\eta M\lambda_j(\Sigma)|<1.
\end{equation}
If $\Sigma$ is singular, coefficients are not uniquely identifiable over this unrestricted feature parameterization.
\end{theorem}
\begin{proof}
By \eqref{eq:square} with $u=k^\T\phi$,
\begin{equation*}
 \mathcal L(k)=\mathcal L_{\rm irr}+M\E[(k^\T\phi-u_\star)^2]
 =\mathcal L_{\rm irr}+M\big(k^\T\Sigma k-2k^\T m_\phi+\E[u_\star^2]\big),
\end{equation*}
a quadratic in $k$ with gradient $\nabla\mathcal L(k)=2M(\Sigma k-m_\phi)$ and Hessian $2M\Sigma$. The gradient vanishes only at $k_\star=\Sigma^{-1}m_\phi$, and $m_\phi=-M^{-1}\E[\phi\,b^\T Ha]$ by the definition of $u_\star$, which gives \eqref{eq:kstar}. Since $\nabla\mathcal L(k_\ell)=2M\Sigma(k_\ell-k_\star)$, a gradient step gives the error recursion $k_{\ell+1}-k_\star=(\Id-2\eta M\Sigma)(k_\ell-k_\star)$. In the eigenbasis of $\Sigma$, the $j$th component of the error is multiplied by $1-2\eta M\lambda_j(\Sigma)$ at every step, and this factor has modulus below one for every $j$ exactly when $0<\eta<1/(M\lambda_{\max}(\Sigma))$, which gives \eqref{eq:rate}. If $\delta\in\ker\Sigma$, then $\E[(\delta^\T\phi)^2]=\delta^\T\Sigma\delta=0$, so $\delta^\T\phi=0$ almost surely, and $k$ and $k+\delta$ define the same action almost surely.
\end{proof}
The same argument holds for empirical averages. With $N$ sampled transitions and a positive-definite empirical Gram matrix $\widehat\Sigma$, the minimizer $\widehat k$ of the empirical cost differs from $k_\star$ only through the sampled features and an average noise term (\cref{app:finite}). When $u_\star$ lies in the PID span and the noise enters as $\xi=bw$ with variance $\sigma^2$, as in experiment E2, this gives the exact mean squared error $\E[N\|\widehat k-k_\star\|^2]=(b^\T Hb/M)^2\sigma^2\operatorname{tr}(\Sigma^{-1})\,N/(N-4)$ for Gaussian resets, because $\E[\widehat\Sigma^{-1}]=N\Sigma^{-1}/(N-4)$.

\Cref{th:recovery} is ordinary least squares: by \eqref{eq:square}, $k_\star$ is the population regression of the one-step optimal action $u_\star$ on the PID features, here induced by a physical cost with no action labels. With the input gain $b$ known, $k_\star$ could also be computed directly by this regression. We analyze gradient descent because it is the procedure that carries over to long horizons, where backpropagation through time replaces the closed form (\cref{sec:long}), and \cref{th:recovery} certifies it where the answer is known. Its convergence is global for the training iteration, from every initial gain vector; it says nothing about the stability of the closed loop that the learned PID produces, which \cref{sec:example} examines.

\Cref{th:real} describes the same controller by its readout weights $\theta=(c_1,c_2,c_0)$. By \eqref{eq:readout}, $c_1=K_i$, $c_2=-K_d/\Delta$, and $c_0=K_p+\Delta K_i+K_d/\Delta$, so $\theta$ is an invertible linear function of $k$, $\theta=\Gamma_\Delta k$, with the matrix below read off row by row. Equivalently, the readout's inputs $(I_{t-1},e_{t-1},e_t)$ are the PID features in other coordinates, $\Gamma_\Delta^{-\T}\phi$, so training the readout solves the same regression. The loss in readout coordinates is $\mathcal L(\Gamma_\Delta^{-1}\theta)$, and by the chain rule its Hessian is $\Gamma_\Delta^{-\T}(2M\Sigma)\Gamma_\Delta^{-1}$. Because the action at every step is linear in $\theta$, convergence $\theta_\ell\to\theta_\star$ makes the network's action converge, at every step and for every error sequence, to that of the PID $k_\star$.
\begin{corollary}[Readout coordinates]\label{cor:readout}
Let $\theta=(c,c_0)=\Gamma_\Delta k$ be the readout parameters of \cref{th:real}, where
\begin{equation}
 \Gamma_\Delta=\begin{bmatrix}0&1&0\\0&0&-1/\Delta\\1&\Delta&1/\Delta\end{bmatrix},\qquad
 \det\Gamma_\Delta=-1/\Delta.
\end{equation}
Under the assumptions of \cref{th:recovery}, the loss in $\theta$ has the unique minimizer $\theta_\star=\Gamma_\Delta k_\star$ and Hessian $\mathcal H_\theta=2M\,\Gamma_\Delta^{-\T}\Sigma\Gamma_\Delta^{-1}$. Gradient descent in $\theta$ converges geometrically from every initialization if and only if $0<\eta<2/\lambda_{\max}(\mathcal H_\theta)$, with rate $q_\theta=\max_j|1-\eta\lambda_j(\mathcal H_\theta)|$. The condition number satisfies
\begin{equation}
 \frac{\kappa(\Gamma_\Delta)^2}{\kappa(\Sigma)}\leq\kappa(\mathcal H_\theta)\leq\kappa(\Sigma)\,\kappa(\Gamma_\Delta)^2,
 \qquad
 \kappa(\Gamma_\Delta)=\frac2\Delta\bigl(1+o(1)\bigr)\quad\text{as }\Delta\to0.
\end{equation}
\end{corollary}
In both coordinates the step bound is $2/\lambda_{\max}$ of the Hessian, $2M\Sigma$ for the gains and $\mathcal H_\theta$ for the readout, so each parameterization needs its own bound. Equivalently, gradient descent in $\theta$ is gradient descent in $k$ preconditioned by $(\Gamma_\Delta^\T\Gamma_\Delta)^{-1}$, so the two follow different trajectories to the same optimum. Because $K_d$ enters the readout scaled by $1/\Delta$, the readout coordinates become ill-conditioned as $\Delta\to0$. At the sampling period $\Delta=0.02$ of the oscillator families of \cref{sec:experiments}, $\kappa(\Gamma_\Delta)\approx100$, so for $\Sigma=\Id_3$ the readout Hessian has condition number about $10^4$ and gradient descent in $\theta$ converges far more slowly along its flattest direction. Its largest eigenvalue is also larger, by the factor $1/\sigma_{\min}(\Gamma_\Delta)^2$, about two at this sampling period, where $\sigma_{\min}$ is the smallest singular value; a step size that is safe for the gains can therefore diverge for the recurrent weights. For these reasons the experiments train the PID part in gain coordinates.

\Cref{app:theory} proves \cref{cor:readout}, plots $\kappa(\Gamma_\Delta)$ over all sampling periods, gives the condition numbers for the example of \cref{sec:example}, and states which gains are recoverable.

\section{Strict improvement while preserving PID coefficients}\label{sec:improvement}
\subsection{Residualization and the exact improvement criterion}
\Cref{th:recovery} finds the best PID, but the one-step optimum $u_\star$ is in general a nonlinear function of the state, and no choice of gains reproduces its nonlinear part. To measure what a larger readout buys, we keep the same memory and add fixed features $\psi(z)\in\R^p$, such as $e^2$:
\begin{equation}\label{eq:extended}
 u_{k,\beta}(z)=k^\T\phi(z)+\beta^\T\psi(z),\qquad\beta\in\R^p.
\end{equation}
The features are fixed because a learned feature extractor would make the problem nonconvex and fall outside the following theorem; the recurrent residual of \cref{sec:warm} learns its features and is tested empirically. The PID submodel is $\beta=0$.

Define the population residualization
\begin{equation}\label{eq:orth}
 C_{\psi\phi}=\E[\psi\phi^\T]\Sigma^{-1},\quad
 \widetilde\psi=\psi-C_{\psi\phi}\phi,\quad
 G=\E[\widetilde\psi\widetilde\psi^\T].
\end{equation}
Here $C_{\psi\phi}$ is the matrix of population regression coefficients of $\psi$ on the PID features: row $j$ of $C_{\psi\phi}\phi$ is the least-squares fit of $\psi_j$ by a linear combination of $e$, $I$, and $D$. The residualized features $\widetilde\psi$ are what that fit leaves, the part of $\psi$ outside the PID span (the set of actions $k^\T\phi$), and $G$ is their Gram matrix. They are orthogonal to the PID features, $\E[\widetilde\psi\phi^\T]=\E[\psi\phi^\T]-C_{\psi\phi}\Sigma=0$. The transformation leaves the joint action class unchanged because $k^\T\phi+\beta^\T\psi=(k+C_{\psi\phi}^\T\beta)^\T\phi+\beta^\T\widetilde\psi$.

\begin{theorem}[Coefficient preservation and strict separation]\label{th:gap}
Assume $\Sigma\succ0$ and $G\succ0$. Parameterize \eqref{eq:extended} using $\widetilde\psi$, and let
\begin{equation}
 m_\psi=\E[\widetilde\psi\,u_\star],\qquad
 \beta_\star=G^{-1}m_\psi.
\end{equation}
The unique optimum is $(k_\star,\beta_\star)$, with exactly the same $k_\star$ as in \cref{th:recovery}. Moreover,
\begin{equation}\label{eq:gap}
 \inf_k\mathcal L(k,0)-\inf_{k,\beta}\mathcal L(k,\beta)
 =Mm_\psi^{\T}G^{-1}m_\psi.
\end{equation}
The gap is strictly positive if and only if $m_\psi\ne0$. The Hessian is $2M\operatorname{diag}(\Sigma,G)$; gradient descent converges geometrically for $0<\eta<1/(M\max\{\lambda_{\max}(\Sigma),\lambda_{\max}(G)\})$.
\end{theorem}
\begin{proof}
With $u=k^\T\phi+\beta^\T\widetilde\psi$, \eqref{eq:square} gives $\mathcal L(k,\beta)=\mathcal L_{\rm irr}+M\E[(k^\T\phi+\beta^\T\widetilde\psi-u_\star)^2]$. When the square is expanded, the cross term $2k^\T\E[\phi\widetilde\psi^\T]\beta$ vanishes by orthogonality, and
\begin{equation*}
 \mathcal L(k,\beta)=\mathcal L_{\rm irr}+M\E[u_\star^2]
 +M\big(k^\T\Sigma k-2k^\T m_\phi\big)
 +M\big(\beta^\T G\beta-2\beta^\T m_\psi\big).
\end{equation*}
The PID block is the quadratic of \cref{th:recovery}, minimized at $k_\star$. The residual block is zero at $\beta=0$ and, completing the square, is minimized at $\beta_\star=G^{-1}m_\psi$, where it equals $-Mm_\psi^\T G^{-1}m_\psi$. The difference of the two minima is \eqref{eq:gap}. The Hessian $2M\operatorname{diag}(\Sigma,G)$ is positive definite, which gives uniqueness and the gradient rate as in \cref{th:recovery}.
\end{proof}
The cross moment $m_\psi$ has a direct reading. Because $\widetilde\psi$ is orthogonal to the PID features,
\begin{equation*}
 m_\psi=\E[\widetilde\psi\,u_\star]=\E\big[\widetilde\psi\,(u_\star-k_\star^\T\phi)\big],
\end{equation*}
and $u_\star-k_\star^\T\phi$ is the part of the optimal action outside the PID span, the residual of the regression of $u_\star$ on $\phi$ that defines $k_\star$. Added capacity therefore helps only when it correlates with the part of the optimal action that the PID span leaves unresolved. Here ``correlates'' means a nonzero cross moment, which is not centered: a constant feature, for example, helps whenever that part has a nonzero mean, as for PID with a bias below. Adding neurons does not by itself imply a strict improvement, and even an unrestricted residual class can improve by at most $M$ times the squared distance of $u_\star$ from the PID span, because by \eqref{eq:square} no policy costs less than $\mathcal L_{\rm irr}$ (\cref{app:theory}). \Cref{th:gap} and \cref{cor:raw} are the Frisch--Waugh--Lovell theorem of partitioned least squares~\citep{frisch1933,lovell1963}, applied to the population regression of $u_\star$ on $(\phi,\psi)$. We parameterize the residual orthogonally so that the PID block keeps its meaning. Its optimum is the one-step optimal PID itself, so switching the residual off ($\beta=0$) leaves that PID, and the learned gains can be compared directly with those of \cref{th:recovery}.

\begin{corollary}[Raw coordinates]\label{cor:raw}
Under the assumptions of \cref{th:gap}, the raw parameterization \eqref{eq:extended} has the unique optimum
\begin{equation}
 \beta_{\rm raw}=\beta_\star,\qquad k_{\rm raw}=k_\star-C_{\psi\phi}^\T\beta_\star.
\end{equation}
Its PID block equals the optimal PID gains $k_\star$ if and only if $C_{\psi\phi}^\T\beta_\star=0$, in particular when $\E[\psi\phi^\T]=0$. If $u_\star=k_{\rm lin}^\T\phi+\beta_{\rm lin}^\T\psi$ exactly, then $k_{\rm raw}=k_{\rm lin}$ for every $\nu$ with a positive-definite joint Gram matrix of $(\phi,\psi)$, equivalently $\Sigma\succ0$ and $G\succ0$, whereas $k_\star=k_{\rm lin}+\Sigma^{-1}\E[\phi\psi^\T]\beta_{\rm lin}$ depends on $\nu$.
\end{corollary}
The corollary follows from the class identity above: the raw pair $(k,\beta)$ produces the same action as the orthogonal pair $(k+C_{\psi\phi}^\T\beta,\beta)$, so the raw optimum is the pair that maps to $(k_\star,\beta_\star)$. Without residualization, the coefficient called ``PID gain'' is therefore $k_\star-C_{\psi\phi}^\T\beta_\star$, and reading that shift as failed recovery would be a mistake. The two parameterizations describe the same controllers; they differ only in which coefficients are called the PID gains, and that assignment depends on $\nu$. In the example below, resets centered at $e=1/2$ make the optimal proportional gain $197/550$, whereas the raw block keeps $K_p=1/2$ (\cref{app:theory}).

\subsection{A nonlinear plant with a closed-form gap}\label{sec:example}
The example is a unit-step discretized double integrator with a quadratic, error-dependent drift. We chose it because a quadratic drift is the simplest nonlinearity that no PID can represent, and because its one-step optimum lies in the span of the PID features and $e^2$, so that every quantity has a closed form. Take $\Delta=1$, constant reference $r_t=0$, and
\begin{equation}\label{eq:example}
 y_{t+1}=2y_t-y_{t-1}+u_t+\alpha e_t^2,\quad \alpha=\frac15,
\end{equation}
with reset states $z=(e,I,d)^\T\sim\mathcal N(0,\Id_3)$, where $d=e_t-e_{t-1}$, and the cost weights $H=\operatorname{diag}(1,1/10,1/10)$ and $R=1$. The weight on $I$ makes the optimal integral gain nonzero; with $H$ acting on $e$ alone it would vanish, because the next error does not depend on $I$.

With $r_t=0$ the error is $e_t=-y_t$, so the plant equation reads $e_{t+1}=2e_t-e_{t-1}-u_t-\alpha e_t^2$. In the coordinates of $z$ the previous error is $e-d$, and one step gives
\begin{equation*}
 e^+=e+d-u-\alpha e^2,\qquad I^+=I+e^+,\qquad d^+=d-u-\alpha e^2,
\end{equation*}
that is, $z^+=a(z)+bu$ with $b=-(1,1,1)^\T$ and a drift $a(z)=Az+\alpha be^2$ that is quadratic in $e$, where $A=\left[\begin{smallmatrix}1&0&1\\1&1&1\\0&0&1\end{smallmatrix}\right]$. The action enters affinely, so (A1) holds, and the transition is noise-free. Here $M=R+b^\T Hb=1+\frac65=\frac{11}{5}$ and $b^\T Ha(z)=-\frac{11}{10}e-\frac1{10}I-\frac65d+\frac65\alpha e^2$, so $u_\star=-M^{-1}b^\T Ha$ from \eqref{eq:square} is the one-step optimum
$u_\star=\frac12e+\frac1{22}I+\frac6{11}d-\frac6{55}e^2$,
whose quadratic term no PID can represent. Take $\phi=z$ and $\psi=e^2$. The moments of the standard Gaussian give every quantity in \cref{th:gap}: $\Sigma=\Id_3$; $\E[\phi e^2]=0$ because odd moments vanish, so $C_{\psi\phi}=0$ and $\widetilde\psi=e^2$; $G=\E[e^4]=3$; and only the quadratic term of $u_\star$ correlates with $e^2$, so $m_\psi=\E[e^2u_\star]=-\frac6{55}\E[e^4]=-\frac{18}{55}$. Hence $k_\star=\Sigma^{-1}m_\phi$ is the linear part of $u_\star$, $\beta_\star=G^{-1}m_\psi$, and the gap between the minimal one-step costs of the PID class, $\mathcal L_{\PID}^\star$, and of the extended class, $\mathcal L_+^\star$, is
\begin{equation}\label{eq:explicit}
 k_\star=\left(\frac12,\frac1{22},\frac6{11}\right)^\T,
 \qquad\beta_\star=-\frac6{55},
\end{equation}
\begin{equation}\label{eq:numericgap}
 \mathcal L_{\PID}^\star-\mathcal L_+^\star
 =\frac{Mm_\psi^2}{G}
 =\frac{11}{5}\,3\left(\frac6{55}\right)^2
 =\frac{108}{1375}=0.07854545\ldots>0.
\end{equation}
These values are derived from the transition and cost; no target gains enter training.

To rule out an improvement due only to a missing action offset, enlarge the baseline to PID with a bias, $u=k^{\T}z+\beta_0$, with minimal cost $\mathcal L_{\PID+\mathrm{bias}}^\star$ and optimal offset $\beta_{0\star}=-6/55$, the mean of $u_\star$. The baseline features are now $(z,1)$. Residualizing the square against them gives $e^2-1$, whose variance is two, because $\E[e^2]=1$ and $\E[ze^2]=0$; its cross moment is $m_\psi=\E[(e^2-1)u_\star]=-\frac{12}{55}$. With $G=2$ this leaves
\begin{equation}\label{eq:biasgap}
 \mathcal L_{\PID+\mathrm{bias}}^\star-\mathcal L_+^\star
 =\frac{Mm_\psi^2}{G}
 =\frac{11}{5}\,2\left(\frac6{55}\right)^2
 =\frac{72}{1375}=0.05236363\ldots>0.
\end{equation}
The optimal PID cost is $\mathcal L_{\PID}^\star=3671/2750$, which is $\mathcal L(0)-M\|k_\star\|^2$ (\cref{app:moments}), so in relative terms the gap is $216/3671\approx5.88\%$ of the optimal PID cost; against PID with bias it is $144/3599\approx4.00\%$. Both gaps grow as $\alpha^2$, because $\beta_\star=-6\alpha/11$; the first is $108\alpha^2/55$, so their size here is set by the choice of $\alpha$. They are also the largest possible on this objective. The extended class contains $u_\star$ exactly, and $u_\star$ minimizes the one-step cost pointwise in the reset state, so no controller of the reset state, however expressive, improves on PID by more than $108/1375$. Residual classes that only approximate $e^2$, or plants whose one-step optimum lies outside the feature span, give smaller gaps.

Neither the gap nor the converged gains transfer unchanged to clipped actuators, bounded reset distributions, or long closed-loop rollouts: the Gaussian has unbounded support, the actuator is unsaturated, and the objective evaluates one transition. In closed loop both the one-step optimal PID and $u_\star$ are locally exponentially stable, because they share the linearization $A+bk_\star^\T$, whose spectral radius is $10/11$ (\cref{prop:closedloop}). Neither is globally asymptotically stable, since both diverge from large enough errors; a quadratic Lyapunov function certifies a region of attraction that holds most Gaussian resets (\cref{fig:lyapunov} in \cref{app:theory}). E3b and E4 (\cref{sec:experiments}) measure closed-loop behavior directly.

\section{From one step to long horizons}\label{sec:long}
\Cref{th:real} concerns the controller alone and holds for every error sequence, in open and closed loop. \Cref{th:recovery,th:gap} instead score one transition from a fixed distribution $\nu$ of the augmented state, in the example $z=(e,I,d)$. In closed loop the visited states depend on the controller, so their distribution is no longer $\nu$, and a one-step optimum is optimal only for the distribution it was fitted on. For the example, \cref{prop:closedloop} gives only local stability. E3b measures this difference on the same plant, and the rest of this section treats the long-horizon objective directly; \cref{tab:experiments} lists which experiment tests which result.

\subsection{Long-horizon objective and practical PID}
For a finite horizon $T$, the long-horizon experiments score a controller by the tracking cost
\begin{equation}\label{eq:longcost}
 J_T(\vartheta)=\E\left[\frac1T\sum_{t=0}^{T-1}
 \left(\widetilde e_t^2+10^{-2}\widetilde{\bar u}_t^2
 +10^{-2}\left(\frac{\bar u_t-\bar u_{t-1}}{s_u}\right)^2\right)
 +\widetilde e_T^2\right],
\end{equation}
where $\vartheta$ collects the controller's trainable parameters, the expectation is over episodes (plant parameters, initial state, reference, noise, and disturbance), $\bar u$ is the applied control, $\widetilde e=e/s_e$ and $\widetilde{\bar u}=\bar u/s_u$ are made dimensionless by fixed scales, $s_e=1$ and $s_u=u_{\max}=2$ (\cref{app:experiments}), $\bar u_{-1}=0$, and the terminal term is not divided by $T$. The squared tracking error dominates. The effort and rate terms, with weight $10^{-2}$, penalize large and rapidly changing actions, and the terminal term penalizes ending away from the reference. Unlike the cost of the example, $J_T$ does not penalize the integral; integral action earns its place by rejecting the actuator offset that E4 applies halfway through each episode. The weights and scales were fixed in the protocol before any training.

Backpropagation through time differentiates this objective exactly, but future states now depend on the parameters, so the fixed-Gram argument of \cref{sec:learning} no longer applies. \Cref{app:longhorizon} gives a conditional convergence result (\cref{prop:longhorizon}) and explains why its assumptions are not established here.

\paragraph{Derivative filtering and anti-windup.}
The practical controllers of \cref{sec:experiments} use
\begin{align}
 D_t&=\lambda_{\rm f}D_{t-1}+(1-\lambda_{\rm f})(e_t-e_{t-1})/\Delta,\\
 v_t&=v_{t-1}+K_i\Delta e_t+k_{\rm aw}(\bar u_{t-1}-u_{t-1}),\\
 u_t&=K_pe_t+v_t+K_dD_t+\rho_t,\qquad
 \bar u_t=\operatorname{clip}(u_t,-u_{\max},u_{\max}).
\end{align}
Here $D_t$ is a filtered derivative that replaces the difference quotient of \eqref{eq:pid}, $v_t$ is the integral already multiplied by $K_i$, $\rho_t$ is the output of the recurrent residual defined in \cref{sec:warm}, and $\rho_t=0$ for a plain PID. The back-calculation term~\citep{astrom1989windup} uses the previous action, which avoids an algebraic loop, and corrects the integral state $v_t$ in control units, which keeps it stable for either sign of $K_i$ (\cref{app:longhorizon}). All compared controllers use the same filtering, clipping, and anti-windup. These updates keep the PID a small recurrent network: the integral $v_t$, the filtered derivative $D_t$, the previous error, and the previous saturation excess form its state, and $K_p$, $K_i$, $K_d$, $\lambda_{\rm f}$ and $k_{\rm aw}$ are its five trainable parameters, of which $K_i$, $\lambda_{\rm f}$ and $k_{\rm aw}$ now enter the memory update. With $\lambda_{\rm f}=k_{\rm aw}=0$, $\rho_t=0$, no clipping, and $v_{-1}=K_iI_{-1}$, it reduces to \eqref{eq:pid}, with $v_t=K_iI_t$. In E4 the realization matters because PID is a submodel of the controller, which makes the residual extension and the warm start below possible, not because of the coordinates in which the gains are trained.

\subsection{Warm start at a tuned PID}\label{sec:warm}
The experiments extend the practical PID with the simplest bounded recurrent network, a single $\tanh$ layer whose output is added to the PID action $u_t$ above:
\begin{equation}\label{eq:residual}
 \rho_t=\rho_{\max}\tanh(W_\rho h_t+b_\rho),\qquad h_t=\tanh(W_hh_{t-1}+W_x\mathbf o_t+b_h).
\end{equation}
Here $\mathbf o_t=(e_t,I_t,D_t,\bar u_{t-1})$ is the controller's observation, with $I_t$ the plain error integral of \eqref{eq:pid} and $D_t$ the filtered derivative. In the experiments $h_t$ has eight units, which gives the controller 118 parameters, the PID's five included. The residual's output can be bounded because the PID part supplies the integral action that, by \cref{prop:lstm}, a bounded network cannot; the amplitude $\rho_{\max}$ sets how far the residual may move the action and is selected on validation. Split the parameters into three groups: the PID parameters $\vartheta_{\rm P}$, namely $K_p$, $K_i$, $K_d$, $\lambda_{\rm f}$ and $k_{\rm aw}$; the memory $\vartheta_{\rm m}=(W_h,W_x,b_h)$; and the output layer $\vartheta_{\rm o}=(W_\rho,b_\rho)$. For one episode, let $p_t$ be the derivative of the episode cost in \eqref{eq:longcost} with respect to an additive perturbation $u_t\leftarrow u_t+\varepsilon_t$ of the unclipped action at $\varepsilon=0$; backpropagation through time computes exactly these sensitivities. Write $J_T^{\PID}(\vartheta_{\rm P})$ for the objective of the practical PID alone.

\begin{proposition}[Warm start]\label{prop:warm}
Let $\vartheta_{{\rm P},0}$ be any PID parameters, $\vartheta_{\rm m}$ arbitrary, and $\vartheta_0=(\vartheta_{{\rm P},0},\vartheta_{\rm m},0)$.
\begin{enumerate}[label=(\roman*)]
\item On every episode, the controller with parameters $\vartheta_0$ produces the same actions and states as the practical PID with $\vartheta_{{\rm P},0}$. Hence $J_T(\vartheta_0)=J_T^{\PID}(\vartheta_{{\rm P},0})$ for every distribution of episodes.
\item If $J_T$ is differentiable at $\vartheta_0$, then $\nabla_{\vartheta_{\rm m}} J_T(\vartheta_0)=0$, $\nabla_{\vartheta_{\rm P}} J_T(\vartheta_0)=\nabla J_T^{\PID}(\vartheta_{{\rm P},0})$, and
\begin{equation}\label{eq:warmgrad}
 \nabla_{W_\rho}J_T(\vartheta_0)=\rho_{\max}\E\Bigl[\sum_tp_th_t\Bigr],\qquad
 \nabla_{b_\rho}J_T(\vartheta_0)=\rho_{\max}\E\Bigl[\sum_tp_t\Bigr],
\end{equation}
with $h_t$ and $p_t$ evaluated on the trajectories of the PID $\vartheta_{{\rm P},0}$.
\item Let $g=\nabla_{\vartheta_{\rm o}}J_T(\vartheta_0)$. If $J_T$ is $L$-smooth, that is, its gradient is $L$-Lipschitz, on a convex set containing $\vartheta_0$ and the point $\vartheta_1$ obtained by the step $\vartheta_{\rm o}\leftarrow-\eta g$, with $0<\eta\leq1/L$, then
$J_T(\vartheta_1)\leq J_T^{\PID}(\vartheta_{{\rm P},0})-\tfrac\eta2\|g\|^2$.
\end{enumerate}
\end{proposition}
The warm start therefore begins exactly at the tuned PID, on every plant and episode, because $\tanh(0)=0$ makes $\rho_t=0$ whatever $h_t$ is. Item (ii) is the chain rule: $\rho_t$ enters the action as the perturbation $\varepsilon_t$ that defines $p_t$, with $\partial\rho_t/\partial W_\rho=\rho_{\max}h_t^\T$, while the path through $h_t$ carries the factor $\rho_{\max}W_\rho=0$. Item (iii) is the standard descent inequality for an $L$-smooth function, applied to a step in the output layer only. A first gradient step therefore improves the training objective whenever $g\neq0$, that is, whenever the sensitivities $p_t$ have a nonzero mean or a nonzero cross moment with the hidden state, and the memory weights receive no gradient until the output layer has moved. The guarantee is first order and concerns the objective whose gradient is computed; it covers neither minibatch Adam over many steps nor generalization to held-out plants, and from an untuned PID it promises improvement only over that start. Differentiability at $\vartheta_0$ is an assumption, because clipping and the rule that aborts diverging rollouts (\cref{sec:experiments}) make $J_T$ piecewise smooth (\cref{app:theory}).

In the one-step setting of \cref{sec:improvement}, the same first step decides whether the residual can improve on the PID. Since $\tanh x=x+O(x^3)$, near $\vartheta_{\rm o}=0$ the residual is $\rho_t\approx\rho_{\max}(W_\rho h_t+b_\rho)$, the extended readout \eqref{eq:extended} with features $\psi=\rho_{\max}(h,1)$ that are not residualized against $\phi$, that is, the raw parameterization of \cref{cor:raw}. The analogy treats $h$ as a fixed feature map of the state, with the memory weights frozen.

\begin{corollary}[One-step warm start]\label{cor:warm}
Under the assumptions of \cref{th:gap}, at the one-step optimal PID $(k,\beta)=(k_\star,0)$, $\nabla_k\mathcal L=0$ and $\nabla_\beta\mathcal L=-2Mm_\psi$, in both the orthogonal and the raw parameterization. A gradient step therefore leaves $k$ unchanged and sets $\beta=2\eta Mm_\psi$, which changes the loss by exactly $-4\eta M^2\|m_\psi\|^2+4\eta^2M^3m_\psi^\T G_\psi m_\psi$, where $G_\psi=G$ in the orthogonal and $G_\psi=\E[\psi\psi^\T]$ in the raw parameterization. For $0<\eta<1/(M\lambda_{\max}(G_\psi))$ the change is negative if and only if $m_\psi\neq0$, that is, if and only if the gap \eqref{eq:gap} is strictly positive.
\end{corollary}
The gradient in $\beta$ at $(k_\star,0)$ is $-2M\E[\psi(u_\star-k_\star^\T\phi)]$, which equals $-2Mm_\psi$ in either parameterization because $u_\star-k_\star^\T\phi$ is orthogonal to $\phi$; the change in loss is then the quadratic $M(\beta^\T G_\psi\beta-2\beta^\T m_\psi)$ at $\beta=2\eta Mm_\psi$. The experiments warm-start at the PID envelope, the best practical PID found by a large search (\cref{sec:experiments}), rather than at $k_\star$, so \cref{cor:warm} is the exact one-step analogue of what \cref{prop:warm} states there to first order; both proofs are in \cref{app:theory}. One difference matters for E4: $k_\star$ is stationary, so the first step leaves the PID part unchanged, whereas the envelope minimizes the long-horizon cost only approximately, so there the PID part moves too (\cref{sec:attribution} measures what that contributes).

\section{Experimental methods}\label{sec:experiments}
Six experiments test the results of \cref{sec:realization,sec:learning,sec:improvement,sec:long} (\cref{tab:experiments}). E1--E3 check the exact results, and E3b connects the one-step theory to closed-loop behavior on the same plant. E5 (E5a for the realization, E5b for the learning) repeats the checks of E1 and E2 in closed loop on classical PID benchmarks; it was added after E4, in response to readers, and like E1--E3 it has pass rules fixed before its reported runs. Their protocols are in \cref{app:experiments}. This section details E4, which asks whether a recurrent residual gives reproducible, worthwhile gains over heavily tuned PID over long horizons, under identical observations and actuator constraints. \Cref{prop:warm} covers only the first step of training, and no theorem covers the trained controllers, so E4 is designed as a confirmatory experiment: its protocol was frozen before any test data existed (\cref{sec:protocol}), and the amendments made after the freeze are numbered 1--4 where they apply.

\begin{table}[tbp]\centering\small
\caption{The experiments, the results they test, the setting of the test (green square: the tested results hold exactly there; amber square: no exact result covers what is measured), and the quantity measured with its pass rule. $N$ is the number of training transitions, and $c_{\rm ex}/N$ the leading term of the learned excess one-step cost (\cref{app:experiments}).}\label{tab:experiments}
\begin{tabularx}{\textwidth}{@{}lll>{\raggedright\arraybackslash}X@{}}
\toprule
Exp. & Tests & Setting & Quantity and rule\\\midrule
E1 & Thm.~\labelcref{th:real}, Prop.~\labelcref{prop:lstm} & \mpass controller only & Paired action discrepancy, at most $10^{-10}$\\\rowgap
E2 & Thm.~\labelcref{th:recovery}, Cor.~\labelcref{cor:readout} & \mpass fixed resets & Gain error against its exact finite-sample law; contraction rates\\\rowgap
E3 & Thm.~\labelcref{th:gap}, Cor.~\labelcref{cor:raw} & \mpass fixed resets & Cost gaps $108/1375$ and $72/1375$; learned excess $c_{\rm ex}/N$\\\rowgap
E3b & Prop.~\labelcref{prop:closedloop} (local) & \mnote closed loop & Closed-loop divergence and cost; no pass rule\\\rowgap
E4 & Prop.~\labelcref{prop:warm} (first step only) & \mnote closed loop & Relative long-horizon cost against the PID envelope\\\rowgap
E5 & Thm.~\labelcref{th:real,th:recovery}, Prop.~\labelcref{prop:warm} & \mpass closed loop, fixed resets & E5a, eight classical plants: discrepancy between PID and network in closed loop, at most $10^{-10}$. E5b, six of them: learned gains against the normal equations, at most $10^{-8}$; trained network against that PID in closed loop, at most $10^{-6}$\\\bottomrule
\end{tabularx}
\end{table}

\subsection{E4: long-horizon closed-loop performance}
\paragraph{Plants and panels.} Three dimensionless benchmark families were used (\cref{tab:families} in \cref{app:e4details}). They differ in how they relate to assumption (A1): the linear oscillator satisfies it, the Duffing oscillator adds a nonlinear drift but keeps the input affine, and the discrete Hammerstein plant passes the input through a saturating nonlinearity, which violates it. The oscillators are simulated by fourth-order Runge--Kutta (RK4) with the control held over $\Delta=0.02$. Episodes last $T=1000$ steps with actuator limit $u_{\max}=2$. From step 500 onward a constant actuator offset $\delta_u$ acts on the plant, which a controller can reject only through integral action, and the observation noise level, and for the oscillators the process noise level, is drawn per episode, so that every controller faces a range of noise (\cref{app:e4details}). Training references mix steps, ramps, pseudo-random binary sequences (PRBS) and multisines with amplitudes up to $0.8$, chosen so that holding them is feasible under $u_{\max}$. Each family has 64 training, 16 validation, and 64 held-out test plants drawn from scrambled Sobol sequences; we call each such set of plants a panel. Every policy is trained across a whole family without access to the plant parameters. Each test plant receives 128 new episodes, 32 per reference type (amendment 1), and all methods see the same episodes. A separate shift panel draws each parameter just outside its training range, with references up to $1.2$ and stronger observation noise; it is reported separately and never pooled with the primary endpoint.

\paragraph{Controllers.} Every controller uses the practical PID of \cref{sec:long}, with derivative filter, anti-windup and clipping, and every learned controller observes $\mathbf o_t=(e_t,I_t,D_t,\bar u_{t-1})$. Three fixed PIDs serve as references. The equal-budget PID is tuned by derivative-free search with about as many simulated episodes as a learned controller receives. The PID envelope is tuned with a far larger budget and approximates the best PID, the lower envelope of the cost over the PID class, hence its name. The PID-plus-bias envelope repeats that search with an added constant offset, which tests whether an offset alone, against the actuator disturbance, explains any gain. Against them we trained the bounded recurrent residual \eqref{eq:residual} with eight units, in two ways that differ only in where and how its PID part starts. The generic start uses parameters drawn once per family from a box of nonnegative values, $K_p,K_i\in[0,2]$, $K_d\in[0,0.5]$, $\lambda_{\rm f}\in[0,0.9]$ and $k_{\rm aw}\in[0,1]$, which have the stabilizing sign because every family has a positive input gain; the warm start uses the PID envelope (\cref{sec:warm}). In both, $W_\rho$ and $b_\rho$ start at zero, all parameters are trained, and the PID part is trained in gain coordinates (\cref{cor:readout}). The PID envelope has 5 parameters, the PID-plus-bias envelope 6, and each residual controller 118.

\paragraph{Training and tuning budgets.} Each learned method was trained with 20 independent seeds, because few-seed comparisons of learned controllers are unreliable~\citep{henderson2018,agarwal2021}. Each run used batches of 16 episodes, a budget of 20,000 simulated episodes, Adam, and gradient clipping at norm 10, and selected among 9 configurations on validation (three learning rates crossed with $\rho_{\max}\in\{0.3,1,2\}$); the selection was recorded with a cryptographic hash before any test evaluation. The residual with a generic start spent the first 20\% of its budget training only its PID part, so that the residual is added to a working PID; this is its own warm start, charged to its budget. The warm-started residual had the full budget, and the compute spent on the envelope is reported separately. The equal-budget PID evaluated 312 Sobol candidates on a 64-episode training panel, 19,968 episodes in total. The PID envelope, by contrast, searched a wide box of gains, filter coefficient, anti-windup fraction and, for the bias variant, offset, combining Sobol sampling~\citep{owen1998}, bounded gradient refinement, and differential evolution~\citep{storn1997} with 10 seeds, and chose its final candidate on validation (\cref{app:e4details}). This intentionally over-budget reference approximates the best PID; it is not a certified global optimum.

\paragraph{Primary endpoint.} For each family $f$, the primary endpoint compares the residual with a generic start against the envelope on the test panel. It asks whether learning, at the learned methods' budget, can replace the over-budget PID tuning; the warm-started residual instead builds on that tuning and is a pre-specified secondary comparison. The frozen protocol made ``worthwhile'' concrete as an improvement of at least 5\%, a judgment of practical relevance fixed before any data, not a value derived from the theory. Two quantities enter: the relative cost $\mathcal R_f=(\overline J_{\rm residual}-\overline J_{\rm envelope})/\overline J_{\rm envelope}$, where $\overline J$ is the mean episode cost over the test plants, their episodes and, for a learned controller, its 20 seeds, and the failure-rate difference $\Delta F_f$ (residual minus envelope), where a failure is an aborted rollout (a state beyond magnitude 20 or a non-finite value, scored with a bounded penalty; \cref{app:e4details}). Failures enter because a lower cost among surviving rollouts would not by itself show that a controller is better. A family is \emph{supported} only if the residual's cost is at least 5\% lower \emph{and} its failure rate at most one percentage point higher. Formally, the one-sided hypotheses $H_{0a}:\mathcal R_f\geq-0.05$ and $H_{0b}:\Delta F_f\geq0.01$ are combined by intersection--union, $p_f=\max(p_a,p_b)$, because a family is supported only if both are rejected, and Holm's procedure~\citep{holm1979} must reject $p_f$ at a familywise error rate of $0.05$ over the three plant families. A family is \emph{not supported} if the one-sided lower bound of $\mathcal R_f$ is at least $-0.05$ or that of $\Delta F_f$ at least $0.01$, and \emph{inconclusive} otherwise (amendment 4). The warm-started residual, the PID-plus-bias envelope, and the equal-budget PID are analyzed with the same estimator, outside Holm.

\subsection{Inference and protocol}\label{sec:protocol}
The $p$-values came from a hierarchical paired bootstrap~\citep{efron1993} with 10,000 replicates. Each replicate resamples test plants, shared by both arms, then episodes within each sampled plant, also shared, and then the training seeds of the learned arm; the envelope is deterministic. The episode is the smallest resampled unit, so time steps are never treated as independent. One-sided percentile bounds were computed at level $1-0.05/3$, a Bonferroni correction over the three families, with bias-corrected and accelerated (BCa) bounds as a sensitivity analysis, and the coverage of the procedure was checked on synthetic data with known $\mathcal R$.

The protocol was frozen on 21 October 2026 (UTC), before any test-panel data existed, and recorded in the repository with a tag; it was not externally preregistered. E1--E3b were run the day before the freeze, under protocol sections that the freeze left unchanged; the freeze binds E4. Four dated amendments and one correction followed the freeze, none using test data (\cref{sec:amend}). Two of them shape what E4 can show. Amendment 3 reduced E4, for compute, to the five controllers above once a validation pilot pointed to a small effect. It removed a gradient-trained PID, generic RNN and LSTM controllers, a quadratic residual, the ablations, and an oracle per-plant PID, so E4 cannot compare the residual with other learned controllers or say which component matters. Amendment 4, which defined the outcome classes \emph{not supported} and \emph{inconclusive}, was committed while the single test evaluation was running, before any of its output was read. One exploratory arm, added after the test evaluation and reported separately, is described in \cref{sec:attribution}.

\section{Results}\label{sec:results}

\subsection{Exact results (E1--E3b, E5)}
Every pass rule of E1--E3 and E5 was met. \Cref{tab:exact} pairs each prediction with its measured value, including the closed-loop behavior that E3b reports without a pass rule; \cref{app:results} gives the full results.

\begin{table}[t]\centering\small
\caption{Exact results: what each experiment predicts and what it measured (green square: pass rule met; amber square: descriptive, no pass rule).}\label{tab:exact}
\begin{tabularx}{\textwidth}{@{}ll>{\raggedright\arraybackslash}X@{}}
\toprule
Exp. & Prediction & Measured\\\midrule
E1 & Realizations of Thm.~\labelcref{th:real} and App.~\labelcref{app:theory} give the same actions & \mpass Agreement to $\EoneMax$ over \EoneCases{} cases\\\rowgap
E1 & A bounded output cannot track the integral (Prop.~\labelcref{prop:lstm}) & \mpass The bound held for every LSTM tested\\\rowgap
E2 & Gradient descent on the physical cost finds $k_\star$ (Thm.~\labelcref{th:recovery}) & \mpass Agrees with the normal equations to $\EtwoNEMaxCorr$\\\rowgap
E2 & Gradient descent in readout coordinates contracts as predicted (Cor.~\labelcref{cor:readout}) & \mpass With exact expectations, rate $\EtwoReadoutRate$, equal to its prediction $1-1/\kappa(\mathcal H_\theta)$\\\rowgap
E2 & Finite-sample gain error follows its exact law (\cref{sec:learning}) & \mpass Mean within $\EtwoMaxAbsZ$ standard errors of the exact value\\\rowgap
E3 & Gaps $108/1375$ and $72/1375$ at $\alpha=1/5$ (Thm.~\labelcref{th:gap}) & \mpass Equal to machine precision: \GapRelPid{} against PID, \GapRelPidBias{} against PID with bias\\\rowgap
E3 & Classes that contain $u_\star$ return $(k_\star,\beta_\star)$ & \mpass To $\EthreeRealizableMax$\\\rowgap
E3b & Local stability only (Prop.~\labelcref{prop:closedloop}) & \mnote The one-step optimal PID diverged in \EthreebDivPid{} of resets, $u_\star$ in none of $10^5$ (\cref{tab:e3b})\\\rowgap
E5a & On classical plants, PID and network give the same closed loop (Thm.~\labelcref{th:real}) & \mpass Agreement to $\EfiveMaxAll$ over \EfiveLoops{} closed loops; the residual at zero is identical to its PID\\\rowgap
E5b & There the network learns the one-step optimal PID (Thm.~\labelcref{th:recovery}) & \mpass Gains to $\EfiveLearnErr$; closed loop to $\EfiveLearnLoop$\\\bottomrule
\end{tabularx}
\end{table}

\begin{figure}[t]
\centering
\includegraphics[width=\textwidth]{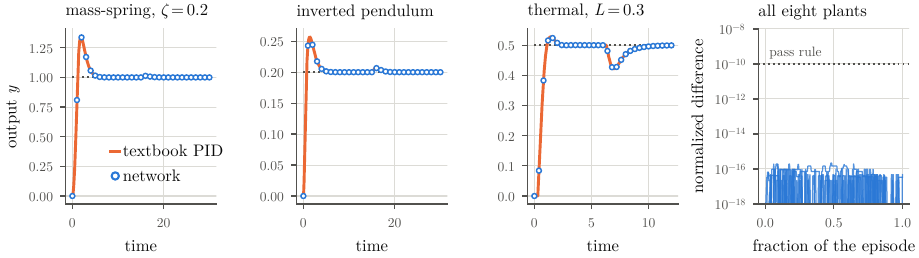}
\caption{E5: closed-loop responses of the textbook PID (line) and of the network of \cref{th:real} (markers) on three of the eight classical plants (all eight in \cref{fig:e5all}), with the classical tuning, actuator saturation, and a load step at mid-episode; the reference is dotted, and time is in the units of each plant model (\cref{tab:e5plants}). The last panel shows their normalized output difference for all eight plants, against the $10^{-10}$ pass rule.}
\label{fig:e5}
\end{figure}

E5 shows the equivalence where a control engineer would look for it, on classical PID benchmarks in closed loop (\cref{fig:e5}). On eight plants, from a mass-spring system, or equivalently a series RLC circuit, in its three damping regimes to an inverted pendulum, two coupled tanks, and a thermal chamber with dead time, and over \EfiveSets{} gain sets and three scenarios each, the network of \cref{th:real} and the textbook PID produce actions and outputs that agree to $\EfiveMaxTwoState$; every realization agrees to $\EfiveMaxAll$, and the E4 controller with its residual at zero is identical to its PID, bit for bit. Because the two loops coincide, so does any quantity computed along them, such as the mechanical energy or a quadratic Lyapunov function of the loop (\cref{fig:e5energy,fig:e5lyap}). What differs is the internal decomposition. The network's readout terms are about an order of magnitude larger than the P, I and D terms and cancel each other (\cref{fig:e5terms}), the ill-conditioning of \cref{cor:readout} seen in the signals. Trained on the one-step physical cost without a teacher, the network reaches the normal-equation gains to $\EfiveLearnErr$ on the six plants where that cost identifies them, and its closed loop then matches that PID's to $\EfiveLearnLoop$. In its readout coordinates it needs at least \EfiveItersNet{} iterations, against at most \EfiveItersGain{} in gain coordinates, both at their predicted rates (\cref{fig:e5conv}); the readout Hessian has condition number about \EfiveCond{}.

\subsection{Long-horizon benchmark (E4)}

\begin{table}[tb]\centering\small
\caption{E4: relative cost $\mathcal R$ against the PID envelope on the test and shift panels (negative is better), with bounds [lower, upper], each one-sided at level $1-0.05/3$, from \EfourBootstrap{} hierarchical paired bootstrap replicates. Green: both bounds below zero; rose: both above zero, or not supported. The decision column (\cref{sec:experiments}) applies to the primary comparison only.}\label{tab:e4}
\begin{tabular}{lrcrcl}
\toprule
 & \multicolumn{2}{c}{warm-started residual} & \multicolumn{2}{c}{residual, generic PID start (primary)} & \\
\cmidrule(lr){2-3}\cmidrule(lr){4-5}
family & $\mathcal R$ (\%) & bounds & $\mathcal R$ (\%) & bounds & decision \\
\midrule
Linear & \cpass $-$3.51 & \cpass [$-$4.51, $-$2.49] & \cpass $-$3.70 & \cpass [$-$4.82, $-$2.49] & \cfail not supported \\
Duffing & \cpass $-$3.45 & \cpass [$-$4.37, $-$2.54] & \cfail $+$14.03 & \cfail [$+$10.92, $+$17.57] & \cfail not supported \\
Hammerstein & \cpass $-$1.90 & \cpass [$-$2.44, $-$1.38] & \cfail $+$3.00 & \cfail [$+$2.21, $+$3.76] & \cfail not supported \\
\midrule
Linear (shift) & \cpass $-$2.01 & \cpass [$-$2.86, $-$1.10] & \cpass $-$2.13 & \cpass [$-$3.18, $-$1.11] & -- \\
Duffing (shift) & \cpass $-$2.30 & \cpass [$-$2.99, $-$1.60] & \cfail $+$17.48 & \cfail [$+$13.51, $+$21.96] & -- \\
Hammerstein (shift) & \cpass $-$2.56 & \cpass [$-$3.49, $-$1.61] & \cfail $+$2.56 & \cfail [$+$1.48, $+$3.71] & -- \\
\bottomrule
\end{tabular}

\end{table}

\paragraph{Reference PIDs.} The PID envelope had a lower validation cost than the equal-budget PID in every family, and on test its relative cost against that PID was \EfourEnvEqLinear{}, \EfourEnvEqDuffing{} and \EfourEnvEqHammerstein{} (linear, Duffing, Hammerstein). Each envelope search used $2\times10^{5}$--$4\times10^{5}$ candidate evaluations of 64 episodes each, about 650 to 1,300 times the budget of one learned training run and 4 to 7 times that of all 180 runs (9 configurations, 20 seeds) of a learned method in the family. The PID-plus-bias envelope scored slightly worse than the PID envelope in every family (\cref{tab:e4sec}), although its class contains every PID. A constant offset therefore does not explain the gains below, and the searches are only approximate. The Hammerstein envelope is the least certain of the three: its anti-windup fraction sits at the physical lower limit $k_{\rm aw}=0$, and its differential-evolution restarts disagree by up to 4.4\% on the training objective.

\paragraph{Primary endpoint: not supported.} No episode failed under any controller, so $\Delta F_f=0$ and the decision rests on $\mathcal R_f$ alone (\cref{tab:e4}; \cref{fig:e4} in \cref{app:results} plots it). In the linear family the residual with a generic PID start lowers the cost ($\mathcal R$ = \EfourBoundedLinear{}), with a point estimate slightly better than that of the warm-started residual, but its one-sided lower bound, \EfourBoundedLoLinear{}, does not reach the $-5\%$ margin. In the Duffing and Hammerstein families it is worse than the envelope ($\mathcal R$ = \EfourBoundedDuffing{} and \EfourBoundedHammerstein{}). Its learning curves (\cref{fig:e4learning}, \cref{app:results}) are consistent with a budget limit: it starts far above the envelope's cost, reaches the warm start's level after about 8,000 episodes in the linear family, and remains above it until the end of training in the other two.

\paragraph{Warm-started residual: a small, consistent gain.} This pre-specified secondary comparison lowers the cost relative to the envelope in all three families, with $\mathcal R$ = \EfourWarmLinear{}, \EfourWarmDuffing{} and \EfourWarmHammerstein{} and one-sided upper bounds of \EfourWarmUpLinear{}, \EfourWarmUpDuffing{} and \EfourWarmUpHammerstein{}. All three bounds exclude zero after a Bonferroni correction over the families, but none reaches the 5\% margin. The improvement persists on the shift panel, with $\mathcal R$ = \EfourWarmShiftLinear{}, \EfourWarmShiftDuffing{} and \EfourWarmShiftHammerstein{}. Shift costs partly measure saturation, however: every test episode is feasible, but 30\% of shift-panel episodes ask for references that $u_{\max}$ cannot hold.

\paragraph{Robustness of the inference.} BCa bounds agree with the percentile bounds to within 0.2 percentage points, and on synthetic data the bootstrap bounds reached their nominal coverage (98.3\% against 98.33\% for $\mathcal R$).

\subsection{What the warm-started residual does}
Little of the gain comes from control effort. The effort term of the cost accounts for only \EfourEffortShareLinear{}, \EfourEffortShareDuffing{} and \EfourEffortShareHammerstein{} of the reduction (linear, Duffing, Hammerstein); the rest comes from the tracking, rate-of-change and terminal terms together, which the evaluation did not record separately. The total variation of the applied action rose in the linear and Hammerstein families (\cref{tab:e4sec}), so the rate term probably did not contribute there. The integrated absolute error (IAE) did not fall. Envelope against warm-started residual, it was about the same in the oscillator families, \EfourIaeEnvLinear{} against \EfourIaeWarmLinear{} and \EfourIaeEnvDuffing{} against \EfourIaeWarmDuffing{}, and higher for Hammerstein, \EfourIaeEnvHammerstein{} against \EfourIaeWarmHammerstein{}. A lower squared cost without a lower IAE is consistent with a residual that reduces large errors at the price of slightly larger small ones, but we did not test this mechanism directly; no controller was tuned for absolute error.

The residual is large. Validation selected $\rho_{\max}=u_{\max}=2$ in both oscillator families, so the residual can span the whole actuator range, and its mean magnitude relative to the mean applied action is \EfourRhoRatioLinear{}, \EfourRhoRatioDuffing{} and \EfourRhoRatioHammerstein{}. It therefore partly opposes the PID command or acts beyond the actuator limit; we did not separate the two cases.

The diagnostics on the training and validation panels confirm \cref{prop:warm}. At initialization the warm-started controller reproduces the envelope's actions exactly and the memory weights receive exactly zero gradient, and along the first gradient direction the change in training cost approaches its first-order prediction, with ratios of observed to predicted change of \EfourFirstOrderLinear{}, \EfourFirstOrderDuffing{} and \EfourFirstOrderHammerstein{} at a step of $10^{-3}$.

\begin{table}[t]\centering\small
\caption{E4, exploratory: the PID tuned alone from the envelope against the envelope, and the warm-started residual against the PID tuned alone, on the test and shift panels (relative cost $\mathcal R$, negative is better, with one-sided bounds at level $1-0.05/3$ from \EfourBootstrap{} hierarchical bootstrap replicates, seeds paired between the two learned arms). The last column is the share of the residual's gain over the envelope that the PID tuned alone reproduces, shown where that PID improves on the envelope. Green: both bounds below zero; rose: both above zero.}\label{tab:attribution}
\begin{tabular}{lrcrcc}
\toprule
 & \multicolumn{2}{c}{PID tuned alone vs envelope} & \multicolumn{2}{c}{residual vs PID tuned alone} & PID's share \\
\cmidrule(lr){2-3}\cmidrule(lr){4-5}
family & $\mathcal R$ (\%) & bounds & $\mathcal R$ (\%) & bounds & of the gain (\%) \\
\midrule
Linear & \cpass $-$1.07 & \cpass [$-$1.60, $-$0.53] & \cpass $-$2.46 & \cpass [$-$3.31, $-$1.60] & 30 [17, 43] \\
Duffing & \cpass $-$0.75 & \cpass [$-$1.10, $-$0.40] & \cpass $-$2.73 & \cpass [$-$3.54, $-$1.95] & 22 [13, 30] \\
Hammerstein & \cpass $-$1.66 & \cpass [$-$2.22, $-$1.08] & \cpass $-$0.25 & \cpass [$-$0.35, $-$0.16] & 87 [78, 93] \\
\midrule
Linear (shift) & \cfail $+$1.59 & \cfail [$+$0.76, $+$2.53] & \cpass $-$3.54 & \cpass [$-$4.44, $-$2.72] & -- \\
Duffing (shift) & \cfail $+$1.28 & \cfail [$+$0.86, $+$1.72] & \cpass $-$3.54 & \cpass [$-$4.30, $-$2.78] & -- \\
Hammerstein (shift) & \cpass $-$2.43 & \cpass [$-$3.47, $-$1.37] & $-$0.13 & [$-$0.34, $+$0.05] & 95 [83, 102] \\
\bottomrule
\end{tabular}

\end{table}

\subsection{Further PID tuning or the recurrent part? (exploratory)}\label{sec:attribution}
The warm-started residual trains its PID parameters along with the residual. Over the 20 seeds, the median change of the proportional and integral gains from the envelope's values stays below 4\%, but that of the filter and anti-windup parameters reaches 29\%, as expected after \cref{cor:warm}: the envelope minimizes the long-horizon cost only approximately, so training moves the PID part too. Part of the gain could therefore come from further tuning of the PID rather than from the recurrent part. To separate the two, we trained the practical PID alone from the same envelope, with the warm-started residual's budget, batch, optimizer, seeds, and training episodes and its learning rate selected on validation from the same three values, and scored it on the same test and shift episodes. The gradient-trained PID that amendment 3 removed would not have answered this, since it started from generic gains. This arm was added after the test panels had been unsealed; its analysis and the rule for reading it were recorded in the repository before it was trained (\cref{app:attribution}). It is therefore reported as post hoc, outside the endpoint.

\Cref{tab:attribution} gives the result. In the two oscillator families the PID tuned alone reaches $\mathcal R$ = \EfourAttrPidLinear{} and \EfourAttrPidDuffing{}, and the residual reaches $\mathcal R$ = \EfourAttrResPidLinear{} and \EfourAttrResPidDuffing{} against it, with upper bounds of \EfourAttrResPidUpLinear{} and \EfourAttrResPidUpDuffing{}; the PID alone reproduces only \EfourAttrShareLinear{} and \EfourAttrShareDuffing{} of the residual's gain. Under shift the PID tuned alone is worse than the envelope in both families, while the residual still improves on it. In the Hammerstein family the picture reverses. The PID tuned alone reproduces \EfourAttrShareHammerstein{} of the gain, and the residual reaches $\mathcal R$ = \EfourAttrResPidHammerstein{} against it on test, with bounds that exclude zero but stay within one percent, and, under shift, a difference that the data do not separate from zero. This is also the family whose envelope is least certain, with restarts that disagree by up to 4.4\%, so its PID left the most room for further tuning. Two limits apply: the comparison is post hoc, and in the linear family the selected learning rate of the PID tuned alone is the largest of the three, so a wider grid might narrow its gap. Within these limits, the data point to the recurrent part, not further tuning of the PID, as the source of most of the gain in the two oscillator families.

\section{Limitations and conclusions}\label{sec:conclusions}
The one-step theory is proved, and E1--E3 confirm it numerically, E2 also statistically; E5 confirms the realization and the recovery in closed loop on classical PID benchmarks. A structured recurrent network contains PID exactly and, trained on a physical cost, learns the one-step optimal gains under full-rank excitation. An orthogonal extension keeps those gains and lowers the cost by a computable amount; on the example plant this proves strict superiority over every fixed PID, with or without bias, for the stated one-step objective.

On held-out plants the evidence is narrower. Started from a heavily tuned PID envelope, a bounded recurrent residual improves on it by 1.9--3.5\% in all three benchmark families and under distribution shift; \cref{prop:warm} shows that it starts exactly at the PID and when its first gradient step lowers the training cost. The protocol's primary claim, an improvement of at least 5\% from a generic PID start, is not supported: that residual matched the warm-started one on the linear family, with a bound that misses the margin, and lost to the envelope on both nonlinear families.

The frozen protocol could not settle attribution, because the warm-started residual trains its PID parameters together with the recurrent part. A post hoc arm that tuned the PID alone from the envelope, at the same budget and on the same episodes, reproduced less than a third of the residual's gain in the two oscillator families, where it did worse than the envelope under shift, but most of the gain in the Hammerstein family (\cref{sec:attribution}). Because it was made after the test panels were unsealed, this comparison supports, but does not confirm, a contribution of the recurrent part in the oscillator families. The gain is also small, and it lowers the squared cost without lowering the integrated absolute error.

The warm-start guarantee is first order and concerns the training objective, and the held-out improvement is empirical. The one-step optimal PID is not globally asymptotically stable on the example plant: a Lyapunov function certifies only a region of attraction, and in E3b it diverged in \EthreebDivPid{} of Gaussian resets. Amendment 3 removed the comparison with other learned controllers and the ablations, the generic start rests on a single random PID per family, and the study covers three simulated families with one sampling period each.

For practice, the results suggest treating a well-tuned PID as the starting point of a recurrent controller rather than as a competitor. At our training budget, starting the PID part at the envelope was needed for any gain on the two nonlinear families, while on the linear family a generic start reached the same cost. In the oscillator families the post hoc comparison suggests that the recurrent part adds to that tuning rather than repeating it; confirming this would need its own frozen protocol.\label{sec:mainend}

\section*{Reproducibility statement}
\addcontentsline{toc}{section}{Reproducibility statement}
\ifpidpreprint The accompanying repository \url{https://github.com/FelixStaff/PIDAsNet}\else An anonymized copy of the public repository, provided as supplementary material,\fi{} contains all code, configurations and analysis scripts. The software environment is locked with \texttt{uv} (Python~3.12, PyTorch with CUDA~12.8); the pinned versions and the reasons for them are documented in the repository.

Every experiment is a configuration file executed by a single command. Each run records its configuration hash, the commit and whether the working tree had uncommitted changes, the software and hardware environment, wall time, work counters, and per-replication or per-episode outputs. Runs from a modified working tree are treated as pilots and never reported. Random streams are derived from documented seeds per experiment, replication and episode. The panels used for selection and those used for the final evaluation are fixed by seeds recorded before any controller is evaluated, and the test panels are generated only after the protocol is frozen.

The analytic constants in this paper are recomputed by the test suite, including symbolic checks of \cref{app:moments} and an independent recomputation of the gaps of \cref{sec:example} that does not use the training code. Training code never imports the module that holds the closed-form optimal gains; a test enforces this. Reported experimental numbers, tables and figures are generated from run outputs by scripts; the few numbers typed in the text are listed with their sources in the repository's evidence map, which a test checks. The experimental protocol, its decisions and any dated amendments are part of the repository.

\section*{Use of AI assistance}
An AI assistant was used to review and proofread the manuscript and to run and verify experiments. The authors directed the work, reviewed the results and the manuscript, and take full responsibility for the content.

\bibliographystyle{tmlr}
\bibliography{refs}

\clearpage
\appendix
\crefalias{section}{appendix}\crefalias{subsection}{subappendix}
\section{Additional theory}\label{app:theory}
This appendix follows the order of \cref{sec:realization,sec:learning,sec:improvement,sec:long}.

\subsection{Other exact realizations}
\paragraph{A current-state realization.}
The dependence on $s_{t-1}$ can be removed by storing one more delayed error. With $\tilde s_t=(I_t,e_t,e_{t-1})^\T$,
\begin{equation}\label{eq:threestate}
 \tilde s_t=\tilde A\tilde s_{t-1}+\tilde b\,e_t,\quad
 \tilde A=\begin{bmatrix}1&0&0\\0&0&0\\0&1&0\end{bmatrix},\quad
 \tilde b=\begin{bmatrix}\Delta\\1\\0\end{bmatrix},\qquad
 u_t=\tilde c^\T\tilde s_t,\quad
 \tilde c=\begin{bmatrix}K_i\\K_p+K_d/\Delta\\-K_d/\Delta\end{bmatrix},
\end{equation}
with $\tilde s_{-1}=(I_{-1},e_{-1},e_{-2})^\T$ for an arbitrary $e_{-2}$, which never reaches the output. Induction gives $\tilde s_t=(I_t,e_t,e_{t-1})^\T$, and expanding $\tilde c^\T\tilde s_t$ yields \eqref{eq:pid}. The map $k\mapsto\tilde c$ has determinant $1/\Delta$, with inverse $K_i=\tilde c_1$, $K_d=-\Delta\tilde c_3$, and $K_p=\tilde c_2+\tilde c_3$. The output is now a linear function of the current state, as in a recurrent network with output $u_t=Wh_t$. The third coordinate is a one-step delay line; in a gated cell it requires a candidate that reads the previously exposed state, because elementwise gates cannot move values between coordinates.

\paragraph{Incremental realization.}
The equivalent recurrence is
\begin{equation}
 u_t=u_{t-1}+a_0e_t+a_1e_{t-1}+a_2e_{t-2},\quad
 \begin{bmatrix}a_0\\a_1\\a_2\end{bmatrix}=
 \begin{bmatrix}1&\Delta&1/\Delta\\-1&0&-2/\Delta\\0&0&1/\Delta\end{bmatrix}k.
\end{equation}
The determinant is one, and
$K_p=-a_1-2a_2$, $K_i=(a_0+a_1+a_2)/\Delta$, $K_d=\Delta a_2$.
Equivalence requires $u_{-1}=K_pe_{-1}+K_iI_{-1}+K_d(e_{-1}-e_{-2})/\Delta$. Here $e_{-2}$ is arbitrary but must be the value used in the first increment, because $a_2e_{-2}$ cancels the $e_{-2}$ term of $u_{-1}$. Any other $u_{-1}$ adds the same constant to every output for a given error sequence. In closed loop, that offset also changes the error sequence, so the trajectories differ by more than a constant.

\paragraph{A gated linear-memory variant.}
Consider the gated cell update $s_t=f\odot s_{t-1}+i\odot g(e_t)$, with identity candidate $g(e)=(\Delta e,e)^\T$ and hard gates $f=(1,0)^\T$ and $i=(1,1)^\T$. It reproduces the two-state recursion of \cref{th:real}, and \eqref{eq:readout} reads the previous cell state $s_{t-1}$ together with the current input $e_t$. This output map differs from the standard LSTM output $h_t=o_t\odot\tanh(c^{\rm cell}_t)$ in two ways: it is linear, and it combines $s_{t-1}$ with a direct input-to-output connection. This is an LSTM-inspired linear-memory variant, not a standard sigmoid/tanh LSTM.

\paragraph{Continuous time.} The ideal continuous-time derivative is not a causal system, so continuous-time PID is used with a first-order derivative filter of time constant $\tau_{\rm f}>0$. That controller is again a fixed linear memory with a three-weight readout.
\begin{proposition}[Continuous-time realization]\label{prop:continuous}
Let $\tau_{\rm f}>0$. With zero initial state, the controller with transfer function $C(s)=K_p+K_i/s+K_ds/(\tau_{\rm f}s+1)$ from $e$ to $u$ is exactly the continuous-time linear recurrent network
\begin{equation}\label{eq:continuous}
 \dot s=A_{\rm c}s+b_{\rm c}e,\quad u=c^\T s+c_0e,\quad
 A_{\rm c}=\begin{bmatrix}0&0\\0&-1/\tau_{\rm f}\end{bmatrix},\
 b_{\rm c}=\begin{bmatrix}1\\1/\tau_{\rm f}\end{bmatrix},\
 c=\begin{bmatrix}K_i\\-K_d/\tau_{\rm f}\end{bmatrix},\
 c_0=K_p+\frac{K_d}{\tau_{\rm f}},
\end{equation}
whose readout is a bijection of the gains, with inverse $K_i=c_1$, $K_d=-\tau_{\rm f}c_2$, and $K_p=c_0+c_2$. If $e$ is held constant over each sampling period $\Delta$, the state at the sampling instants obeys $s_{t+1}=A_{\rm d}s_t+b_{\rm d}e_t$, with $A_{\rm d}=\operatorname{diag}(1,\gamma)$, $b_{\rm d}=(\Delta,1-\gamma)^\T$ and $\gamma=\exp(-\Delta/\tau_{\rm f})$, and the action is $u_t=c^\T s_t+c_0e_t$ with the same readout.
\end{proposition}
\begin{proof}
The first state is the integral of $e$, and the second is $e$ filtered by $1/(\tau_{\rm f}s+1)$. Then $c^\T s+c_0e=K_pe+K_is_1+K_d(e-s_2)/\tau_{\rm f}$, and $(e-s_2)/\tau_{\rm f}$ has transfer function $\big(1-1/(\tau_{\rm f}s+1)\big)/\tau_{\rm f}=s/(\tau_{\rm f}s+1)$ from $e$, which gives $C(s)$. The map from $(K_p,K_i,K_d)$ to $(c,c_0)$ is linear with the stated inverse. Over a period with constant $e_t$, $\dot s_1=e_t$ adds $\Delta e_t$, and $\dot s_2=(e_t-s_2)/\tau_{\rm f}$ moves $s_2$ to $\gamma s_2+(1-\gamma)e_t$.
\end{proof}
The readout has the structure of \cref{th:real}, with the filter time constant in the role of the sampling period: the derivative gain enters divided by $\tau_{\rm f}$, as it enters divided by $\Delta$ there. A test of the repository checks the transfer function, the inverse, and the sampled recursion symbolically.

\subsection{Bounded-output readouts and nested optima}
\Cref{prop:lstm} has a narrow scope. For $K_i=0$ and bounded errors, a PD output is bounded, so the proposition does not exclude exact realization in that case. The argument concerns a standard exposed hidden-state readout; direct access to an unbounded cell state changes the architecture. It also concerns unsaturated exact realization only. Under the actuator clip $|\bar u_t|\leq u_{\max}$ of \cref{sec:long}, the applied action of every controller is bounded, so the proposition places no restriction on what an LSTM can apply there. The practically relevant question is whether its cell state learns to integrate, which is empirical. Finite-horizon approximation is compatible with the proposition, but neither exact integration nor training success follows from approximation alone.

Exact realization also makes PID and its recurrent realization interchangeable as policies. Let $J(\pi)$ be the expected trajectory cost of a policy $\pi$, $\pi_k$ the PID with gains $k$, $\pi_\theta$ the recurrent realization with readout parameters $\theta\in\Theta_{\PID\text{-RNN}}$, and $\Pi_+$ a larger policy class.
\begin{proposition}[Equivalence and nested performance]\label{prop:nesting}
For the same plant, initial conditions, observations, actuator, and trajectory cost, the exact PID and its recurrent realization have identical risks. If $\Pi_+$ contains every PID as a zero-residual special case, then
\begin{equation}
 \inf_{\pi\in\Pi_+}J(\pi)\leq\inf_{k\in\R^3}J(\pi_k)
 =\inf_{\theta\in\Theta_{\PID\text{-RNN}}}J(\pi_\theta).
\end{equation}
\end{proposition}
\begin{proof}
Couple exogenous randomness. Equal actions and equal states imply equal next states by induction, hence equal trajectory costs. Taking expectations and using set inclusion proves the claims.
\end{proof}
This is a statement about attainable optima, not about the outcome of a finite training run.

\subsection{Readout coordinates: proof of \texorpdfstring{\cref{cor:readout}}{the readout-coordinates corollary} and conditioning}
\begin{figure}[t]
\centering
\includegraphics[width=0.62\textwidth]{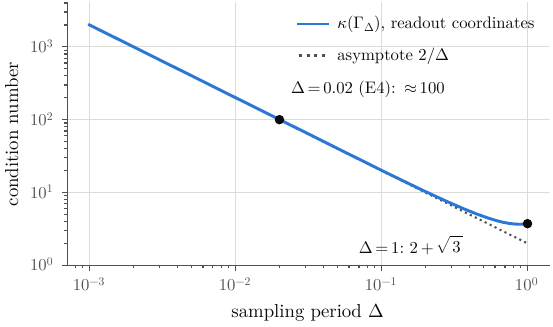}
\caption{Condition number of the readout map $\Gamma_\Delta$ (solid) and its asymptote $2/\Delta$ (dotted). With $\Sigma=\Id_3$ the Hessian in readout coordinates has condition number $\kappa(\Gamma_\Delta)^2$ (\cref{cor:readout}). Computed numerically from $\Gamma_\Delta$; no experimental data.}
\label{fig:cond}
\end{figure}

\begin{proof}[Proof of \cref{cor:readout}]
The map $\theta\mapsto\Gamma_\Delta^{-1}\theta$ is a linear bijection, so the loss in $\theta$ is the quadratic of \cref{th:recovery} composed with it. This gives $\theta_\star$ and $\mathcal H_\theta$, and the recursion $\theta_{\ell+1}-\theta_\star=(\Id-\eta\mathcal H_\theta)(\theta_\ell-\theta_\star)$ gives the rate. The bounds follow from $\lambda_{\min}(\Sigma)\,\Gamma_\Delta^{-\T}\Gamma_\Delta^{-1}\preceq\Gamma_\Delta^{-\T}\Sigma\Gamma_\Delta^{-1}\preceq\lambda_{\max}(\Sigma)\,\Gamma_\Delta^{-\T}\Gamma_\Delta^{-1}$ and $\kappa(\Gamma_\Delta^{-1})=\kappa(\Gamma_\Delta)$. Multiplied by $\Delta^2$, the characteristic polynomial of $\Gamma_\Delta\Gamma_\Delta^\T$ is $\Delta^2\lambda^3-(\Delta^4+2\Delta^2+2)\lambda^2+(2\Delta^2+3)\lambda-1$. Two of its roots tend to $1/2$ and $1$ as $\Delta\to0$. The product of all three roots is $\det(\Gamma_\Delta\Gamma_\Delta^\T)=\Delta^{-2}$, so the largest root is $2\Delta^{-2}(1+o(1))$.
\end{proof}
In the example of \cref{sec:example} ($\Delta=1$, $\Sigma=\Id_3$), the eigenvalues of $\Gamma_1\Gamma_1^\T$ are $2-\sqrt3$, $1$, and $2+\sqrt3$. Hence $\kappa(\mathcal H_\theta)=7+4\sqrt3\approx13.93$ against $\kappa=1$ in gain coordinates, and $\lambda_{\max}(\mathcal H_\theta)=\tfrac{22}{5}(2+\sqrt3)\approx16.42$ against $22/5$. \Cref{fig:cond} plots $\kappa(\Gamma_\Delta)$ over all sampling periods.

\subsection{Scope of the recovery theorem}
\paragraph{Which gains are recovered.}
The recoverable vector is the optimal PID vector for the declared plant, distribution, and physical objective. It is not an arbitrary engineer-selected PID, nor a universal coefficient vector independent of the cost weights $H$ and $R$. Without observing an arbitrary controller or specifying an objective that identifies it, recovering its particular gains is impossible: the same training observations could be paired with any unobserved gain label.

\paragraph{State-dependent input gain and gradients without a model.}
A state-dependent $b(z)$ in (A1) preserves the argument of \cref{sec:learning} with the weight $M(z)=R+b(z)^\T Hb(z)$ inside the expectation in \eqref{eq:square}, which turns the projections into weighted ones. Without a differentiable transition model, the gradient of the quadratic objective \eqref{eq:physical} can still be obtained exactly by central differences in $k$, using costs of alternative actions at the same reset states and the same noise realizations. Without such common random numbers, derivative-free estimates add variance, and a learned simulator adds model error.

\paragraph{Linear augmented dynamics.}
If $a(z)=Az$ and $\phi=z\in\R^3$, then
\begin{equation}\label{eq:linearopt}
 k_\star^\T=-M^{-1}b^{\T}HA,
\end{equation}
independently of the positive-definite design covariance. This is a one-step quadratic controller expressed in PID coordinates; it is not automatically the infinite-horizon linear quadratic regulator (LQR) solution.

\subsection{Proof of \texorpdfstring{\cref{cor:raw}}{the raw-coordinates corollary} and the unrestricted residual limit}
\begin{proof}[Proof of \cref{cor:raw}]
By the class identity after \eqref{eq:orth}, the raw coefficients $(k,\beta)$ and the orthogonal coefficients $(k+C_{\psi\phi}^\T\beta,\beta)$ define the same action, and the map between them is a linear bijection. The unique orthogonal optimum $(k_\star,\beta_\star)$ therefore corresponds to the stated raw optimum. In the realizable case, $\E[\phi u_\star]=\Sigma k_{\rm lin}+\E[\phi\psi^\T]\beta$, so $k_\star=k_{\rm lin}+C_{\psi\phi}^\T\beta$. Moreover $m_\psi=\E[\widetilde\psi u_\star]=G\beta$ because $\E[\widetilde\psi\phi^\T]=0$, so $\beta_\star=\beta$.
\end{proof}
\paragraph{Unrestricted residual limit.}
In the Hilbert space $L^2(\nu)$, let $P_\phi$ be the orthogonal projection onto the PID span. If the residual class contains $u_\star-P_\phi u_\star$, then the attainable improvement is exactly $M\|u_\star-P_\phi u_\star\|_{L^2(\nu)}^2$. A finite feature class recovers the projection of this discrepancy onto its own span. No residual class does better: by \eqref{eq:square} every controller of $z$ has cost at least $\mathcal L_{\rm irr}$, whereas the best PID has cost $\mathcal L_{\rm irr}+M\|u_\star-P_\phi u_\star\|_{L^2(\nu)}^2$.

\subsection{The example in detail}
At each design sample, let $z=(e,I,d)^\T\sim\mathcal N(0,\Id_3)$, where $d=e_t-e_{t-1}$ and $I=I_t$. Such reset states are consistent: set $y_t=-e$, $y_{t-1}=d-e$, and $I_{t-1}=I-e$. These are independent resets, not claims about a Gaussian on-policy stationary distribution. Updating the actual features gives
\begin{equation}
 e^+=e+d-u-\alpha e^2,\qquad I^+=I+e^+,\qquad d^+=d-u-\alpha e^2.
\end{equation}
Consequently
\begin{equation}
 z^+=Az+b(u+\alpha e^2),\quad
 A=\begin{bmatrix}1&0&1\\1&1&1\\0&0&1\end{bmatrix},\quad
 b=\begin{bmatrix}-1\\-1\\-1\end{bmatrix}.
\end{equation}
Choose $H=\operatorname{diag}(1,1/10,1/10)$ and $R=1$. Then
\begin{equation}
 M=\frac{11}{5},\quad b^{\T}HA=(-11/10,-1/10,-6/5),\quad
 u_\star=\frac12e+\frac1{22}I+\frac6{11}d-\frac6{55}e^2.
\end{equation}
For $\phi=z$ and $\psi=e^2$, centered Gaussian symmetry gives
$\E[\phi e^2]=0$, $\Sigma=\Id_3$, and $G=\E[e^4]=3$. Substituting these moments into \cref{th:gap} gives \eqref{eq:explicit} and \eqref{eq:numericgap}.

\begin{figure}[t]
\centering
\includegraphics[width=\textwidth]{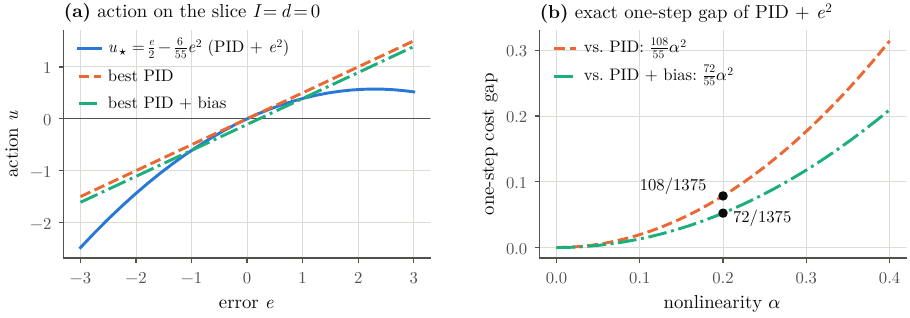}
\caption{The example of \cref{sec:example}. (a)~On the slice $I=d=0$, the one-step optimum $u_\star$ bends away from the best PID and the best PID with bias; the difference is the quadratic term that no PID can represent. (b)~Exact one-step cost gaps of PID plus $e^2$ against PID and against PID with bias, as functions of $\alpha$; the markers are \eqref{eq:numericgap} and \eqref{eq:biasgap} at $\alpha=\frac15$. Analytic; no experimental data.}
\label{fig:example}
\end{figure}

One extra square-activation feature is sufficient. A ReLU or tanh approximation to the square must be evaluated as an approximation and does not inherit exact polynomial equality. \Cref{fig:example} shows where the gap comes from and how it scales with $\alpha$.

Without residualization, the raw PID block shifts with the design distribution (\cref{cor:raw}). Centered Gaussian resets give $C_{\psi\phi}=0$, so here the two parameterizations agree. With resets $z\sim\mathcal N((1/2,0,0)^\T,\Id_3)$, the optimal PID proportional gain becomes $197/550$, whereas the raw PID-plus-$e^2$ controller keeps $K_p=1/2$, the linear part of $u_\star$.

The one-step optima are fitted to independent resets; the following proposition describes their closed loop near the origin.
\begin{proposition}[Closed loop of the one-step optimum]\label{prop:closedloop}
For $\alpha=0$, the closed loop under the one-step optimal PID gains \eqref{eq:explicit} is $z_{t+1}=(A+bk_\star^\T)z_t$, with characteristic polynomial $(\lambda-\tfrac{10}{11})(\lambda^2-\lambda+\tfrac12)$. Its eigenvalues are $\tfrac12\pm\tfrac12i$ and $\tfrac{10}{11}$, so it is asymptotically stable with spectral radius $10/11$. For $\alpha>0$, the origin is an equilibrium of the closed loops under both $k_\star$ and $u_\star$, and both are locally exponentially stable with this same linearization.
\end{proposition}
\begin{proof}
Direct computation gives
\begin{equation*}
 A+bk_\star^\T=\begin{bmatrix}1/2&-1/22&5/11\\1/2&21/22&5/11\\-1/2&-1/22&5/11\end{bmatrix}
\end{equation*}
and the stated characteristic polynomial, and $|\tfrac12\pm\tfrac12i|=1/\sqrt2<10/11$. For $\alpha>0$, the closed-loop maps are $z\mapsto(A+bk_\star^\T)z+b(\alpha+\beta)e^2$, with $\beta=0$ under $k_\star$ and $\beta=\beta_\star$ under $u_\star$, so the quadratic coefficient is $\alpha$ or $\alpha+\beta_\star=5\alpha/11$. Both vanish at $z=0$ and have Jacobian $A+bk_\star^\T$ there, and a fixed point whose Jacobian has spectral radius below one is locally exponentially stable.
\end{proof}
\begin{remark}
The one-step optimum is stabilizing in this example, but it is not the infinite-horizon optimum. The discrete-time LQR for the same $(A,b,H,R)$, computed numerically from the Riccati equation, has gains $(0.558,0.130,0.831)$ and closed-loop spectral radius $0.730$. Stability of one-step optima is not guaranteed in general: for $x_{t+1}=2x_t+u_t$ with a scalar weight $H<1$ and $R=1$, the one-step optimum gives the closed-loop factor $2/(1+H)>1$. For $\alpha>0$ the paragraph below bounds the regions of attraction.
\end{remark}

\begin{figure}[t]
\centering
\includegraphics[width=\textwidth]{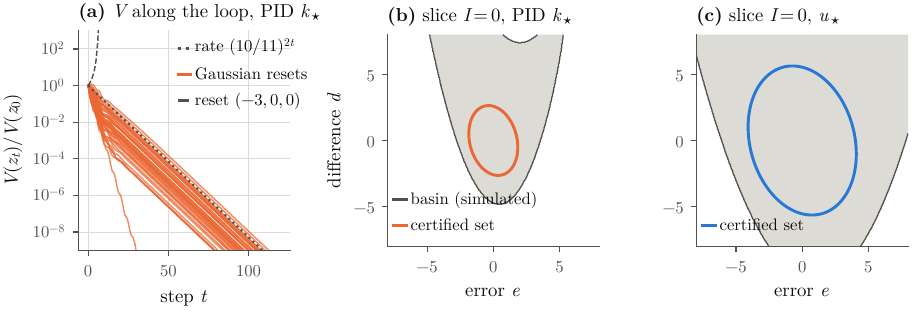}
\caption{Lyapunov analysis of the example closed loop. (a)~$V(z_t)/V(z_0)$ under the one-step optimal PID from 40 Gaussian resets (orange), approaching the rate $(10/11)^{2t}$ of the linearization (dotted), and from one reset outside the basin of attraction (dashed). (b, c)~Slice $I=0$ of the state space: the basin of attraction over 200 steps, computed by simulation (gray), and the set certified by $V$ (line), under $k_\star$ and under $u_\star$. Computed from the model; no experimental data.}
\label{fig:lyapunov}
\end{figure}

\paragraph{Lyapunov analysis and regions of attraction.} \Cref{prop:closedloop} is local, and neither closed loop is globally asymptotically stable: the PID $k_\star$ diverges from the state $(e,I,d)=(\LyapDivPid,0,0)$, and $u_\star$ from $(\LyapDivStar,0,0)$, because the quadratic drift then outgrows the linear feedback. A quadratic Lyapunov function certifies a region of attraction around the origin. Write $A_{\rm cl}=A+bk_\star^\T$ and $c=b(\alpha+\beta)$, so that the closed loop is $z^+=A_{\rm cl}z+ce^2$. For a weight $Q\succ0$, let $P\succ0$ solve the discrete Lyapunov equation $A_{\rm cl}^\T PA_{\rm cl}-P=-Q$, which has a solution because $A_{\rm cl}$ is stable, and let $V(z)=z^\T Pz$. Along the closed loop,
\begin{equation*}
 V(z^+)-V(z)=-z^\T Qz+2e^2\,c^\T PA_{\rm cl}z+e^4\,c^\T Pc.
\end{equation*}
On a ray $z=r\omega$ with $\|\omega\|=1$ this equals $r^2(-\omega^\T Q\omega+2a_\omega r+b_\omega r^2)$, with $a_\omega=\omega_e^2c^\T PA_{\rm cl}\omega$ and $b_\omega=\omega_e^4c^\T Pc\geq0$, where $\omega_e$ is the error component of $\omega$, so it is negative exactly for $0<r<r(\omega)$, the positive root of the quadratic. The largest sublevel set $\{V<\varrho_\star\}$ on which $V$ decreases therefore has $\varrho_\star=\min_\omega r(\omega)^2\omega^\T P\omega$. It is invariant, and every trajectory that starts in it converges to the origin. We compute the minimum over a dense set of directions with local refinement and choose a diagonal $Q$ that maximizes the Gaussian mass of the set; any $Q\succ0$ gives a valid certificate. The certified set contains \LyapMassPid{} of Gaussian resets under $k_\star$ and \LyapMassStar{} under $u_\star$, whereas in E3b only \EthreebDivPid{} of resets diverge under $k_\star$ and none under $u_\star$ (\cref{tab:e3b}); the difference is the conservatism of a quadratic certificate. \Cref{fig:lyapunov} shows $V$ decreasing along the closed loop at the rate of the linearization, and the certified sets inside the basins of attraction computed by simulation.

\subsection{Proofs for \texorpdfstring{\cref{sec:warm}}{the warm-start section}}
\begin{proof}[Proof of \cref{prop:warm}]
(i) With $\vartheta_{\rm o}=0$, $\rho_t=\rho_{\max}\tanh(0)=0$ for every $h_t$, and the remaining recursions are those of the practical PID. (ii) At $\vartheta_{\rm o}=0$ the residual does not depend on $\vartheta_{\rm m}$, and by (i) neither does the trajectory, so $\nabla_{\vartheta_{\rm m}} J_T=0$; the trajectory and its $\vartheta_{\rm P}$-derivative are those of the PID. The residual enters as the perturbation $\varepsilon_t=\rho_t$, with $\partial\rho_t/\partial W_\rho=\rho_{\max}h_t^\T$ and $\partial\rho_t/\partial b_\rho=\rho_{\max}$ at $\vartheta_{\rm o}=0$. The indirect dependence of $\rho_t$ on $\vartheta_{\rm o}$ through $h_t$ carries the factor $\partial\rho_t/\partial h_t=\rho_{\max}W_\rho=0$, and every downstream effect of $\rho_t$ is contained in $p_t$. The chain rule, with differentiation under the expectation, gives \eqref{eq:warmgrad}. (iii) This is the descent inequality used in \cref{prop:longhorizon}, applied to a step in $\vartheta_{\rm o}$ only.
\end{proof}
Clipping and the abort rule make $J_T$ piecewise smooth, so differentiability at $\vartheta_0$ is an assumption of \cref{prop:warm}; it holds when, with probability one, no action lies exactly at $\pm u_{\max}$ and no state lies exactly at the abort threshold.
\begin{proof}[Proof of \cref{cor:warm}]
By \eqref{eq:square}, $\nabla_\beta\mathcal L(k_\star,0)=2M\E[\psi(k_\star^\T\phi-u_\star)]$. The normal equations give $\E[\phi(k_\star^\T\phi-u_\star)]=0$, so $\psi$ may be replaced by $\widetilde\psi=\psi-C_{\psi\phi}\phi$. Since $\E[\widetilde\psi\phi^\T]=0$, the gradient equals $-2M\E[\widetilde\psi u_\star]=-2Mm_\psi$. Since $\nabla_k\mathcal L=0$, the step changes only $\beta$. The loss is quadratic in $\beta$: $\mathcal L(k_\star,\beta)-\mathcal L(k_\star,0)=M(\beta^\T G_\psi\beta-2\beta^\T m_\psi)$, where the quadratic term is $\E[(\beta^\T\widetilde\psi)^2]=\beta^\T G\beta$ in the orthogonal and $\E[(\beta^\T\psi)^2]=\beta^\T\E[\psi\psi^\T]\beta$ in the raw parameterization, and the linear term is the gradient above in both. Substituting $\beta=2\eta Mm_\psi$ gives the stated change, and $m_\psi^\T G_\psi m_\psi\leq\lambda_{\max}(G_\psi)\|m_\psi\|^2$ gives its sign.
\end{proof}

\section{Exact moments and independent checks}\label{app:moments}
For the Gaussian example, $\E[e^2]=1$, $\E[e^3]=0$, $\E[e^4]=3$, and mixed odd moments vanish. The joint feature vector $(\phi^\T,\psi)^\T=(e,I,d,e^2)^\T$ has Gram matrix $\operatorname{diag}(1,1,1,3)$, hence joint Hessian
\begin{equation}
 \nabla^2\mathcal L=\operatorname{diag}(22/5,22/5,22/5,66/5).
\end{equation}
Full-batch population gradient descent converges for $0<\eta<5/33$. With $\eta=5/66$, the residual coefficient reaches its optimum in one update, and the remaining coefficient errors contract by $2/3$ per step. This statement uses exact population gradients.

The zero-action cost is
\begin{equation}
 \E[(Az+\alpha b e^2)^{\T}H(Az+\alpha b e^2)]
 =\operatorname{tr}(A^{\T}HA)+\alpha^2(b^{\T}Hb)\E[e^4]
 =\frac{318}{125}.
\end{equation}
The cross term vanishes. Since $\|k_\star\|^2=133/242$, the exact optimal costs are
\begin{align}
 \mathcal L_{\PID}^\star&=\frac{318}{125}-\frac{133}{110}
 =\frac{3671}{2750},\\
 \mathcal L_+^\star&=\frac{3671}{2750}-\frac{108}{1375}
 =\frac{691}{550},\\
 \mathcal L_{\PID+\mathrm{bias}}^\star&=\frac{3671}{2750}-\frac{36}{1375}
 =\frac{3599}{2750}.
\end{align}
These constants permit a direct independent implementation check against the physical transition cost. They are reported as analytic values, not as estimated experimental means.

\section{Finite-sample coefficient perturbation}\label{app:finite}
Let $m_\phi=\E[\phi u_\star]$, $\widehat\Sigma=N^{-1}\sum_i\phi_i\phi_i^\T$, $\widehat m_\phi=N^{-1}\sum_i\phi_i u_\star(z_i)$, and
\begin{equation}
 \widehat\epsilon=-M^{-1}N^{-1}\textstyle\sum_i\phi_i\,b^\T H\xi_i,
\end{equation}
where $\xi_i$ is the noise realized in the $i$-th observed transition. These quantities are analysis notation; training minimizes the empirical physical cost directly. Its gradient is $2N^{-1}\sum_i\phi_i[Mk^\T\phi_i+b^\T H(a(z_i)+\xi_i)]$, so when $\widehat\Sigma\succ0$ its unique minimizer $\widehat k$ satisfies $\widehat\Sigma\widehat k=\widehat m_\phi+\widehat\epsilon$. If $\lambda_{\min}(\Sigma)\geq\mu>0$ and $\|\widehat\Sigma-\Sigma\|_2\leq\mu/2$, then
\begin{equation}
 \|\widehat k-k_\star\|_2\leq\frac2\mu
 \left(\|\widehat m_\phi-m_\phi\|_2+\|\widehat\epsilon\|_2+\|\widehat\Sigma-\Sigma\|_2\|k_\star\|_2\right).
\end{equation}
Indeed $\widehat\Sigma(\widehat k-k_\star)=(\widehat m_\phi-m_\phi)+\widehat\epsilon-(\widehat\Sigma-\Sigma)k_\star$, and Weyl's inequality gives $\lambda_{\min}(\widehat\Sigma)\geq\mu/2$, so the inverse norm is at most $2/\mu$. For independent samples, (A2) gives $\E[\widehat\epsilon]=0$ and $\operatorname{Cov}(\widehat\epsilon)=M^{-2}N^{-1}\E[(b^\T H\xi)^2\phi\phi^\T]$; the noise term vanishes for noiseless transitions. If, in addition, $u_\star$ lies in the span of $\phi$, as in the linear augmented case \eqref{eq:linearopt}, then $\widehat m_\phi=\widehat\Sigma k_\star$ and $\widehat k-k_\star=\widehat\Sigma^{-1}\widehat\epsilon$ exactly. Without noise the empirical minimizer then equals $k_\star$ for every $N$ with $\widehat\Sigma\succ0$. This deterministic bound separates statistical moment error from optimization error. No distribution-free finite-sample rate is asserted without appropriate tail assumptions.

\section{Long-horizon convergence and anti-windup}\label{app:longhorizon}
This appendix states the conditional convergence result referred to in \cref{sec:long} and explains the anti-windup law shared by every practical controller.

\paragraph{Conditional convergence.}
Let $\chi_{t+1}=F(\chi_t,\vartheta,\xi_t)$ collect plant and controller states. Exact forward sensitivities satisfy
\begin{equation}
 S_{t+1}=F_\chi S_t+F_\vartheta,\qquad S_0=0.
\end{equation}
Backpropagation through time evaluates the same chain rule in reverse.

\begin{proposition}[Conditional long-horizon convergence]\label{prop:longhorizon}
Let $\mathcal V$ be a region containing all gradient-descent iterates and the intervening update segments, and let $J_T^\star=\inf_{\vartheta\in\mathcal V}J_T(\vartheta)$. Suppose $J_T$ is $L$-smooth on $\mathcal V$ and satisfies the Polyak--\L{}ojasiewicz inequality
$\frac12\|\nabla J_T\|^2\geq\mu(J_T-J_T^\star)$ there. For $0<\eta\leq1/L$,
\begin{equation}
 J_T(\vartheta_\ell)-J_T^\star\leq(1-\eta\mu)^\ell
 [J_T(\vartheta_0)-J_T^\star].
\end{equation}
If also $J_T(\vartheta)-J_T^\star\geq(\mu_{\rm qg}/2)\|\vartheta-\vartheta_\star\|^2$ for a unique minimizer, the squared parameter error is bounded by $2/\mu_{\rm qg}$ times this objective bound.
\end{proposition}
\begin{proof}
Smoothness gives $J_T(\vartheta-\eta\nabla J_T)\leq J_T(\vartheta)-(\eta/2)\|\nabla J_T\|^2$. Apply the stated gradient inequality recursively, then the quadratic-growth assumption.
\end{proof}
This is the standard Polyak--\L{}ojasiewicz argument~\citep{polyak1963,karimi2016}. With $J_T^\star$ defined over $\mathcal V$, the bound certifies convergence only to the best value in that region; it is a global statement only if $\mathcal V$ contains a global minimizer. The substance lies in verifying the assumptions. For static state feedback on linear plants with quadratic cost, the infinite-horizon LQR cost satisfies a gradient-dominance condition over the set of stabilizing gains, which yields global convergence of policy gradient~\citep{fazel2018}. That structure does not transfer automatically to restricted controller classes. For static decentralized controllers, the stabilizing set can have exponentially many connected components, so local search from an arbitrary initialization need not reach a global optimum~\citep{feng2020}. PID policies are structured dynamic output feedback; \citet{sharma} establish guarantees in their setting. The assumptions of the proposition are not established for arbitrary nonlinear plants, generic LSTMs, or the benchmark families of \cref{sec:experiments}. Numerical Hessians near a solution provide local diagnostics, not a global proof.

\begin{figure}[t]
\centering
\includegraphics[width=0.62\textwidth]{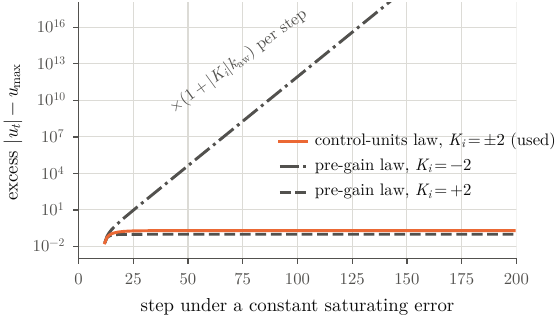}
\caption{Saturation excess $|u_t|-u_{\max}$ under a constant saturating error for $K_i=\pm2$, plotted while the actuator saturates. Correcting the pre-gain integral works for $K_i>0$ but diverges geometrically for $K_i<0$; the control-units law gives the same bounded excess for both signs. Deterministic simulation of both laws with the parameters of the test suite; no experimental data.}
\label{fig:antiwindup}
\end{figure}

\paragraph{Anti-windup in control units.} For the practical controller of \cref{sec:long}, the integral state $v_t$ is kept in control units, so the correction changes $u$ by $k_{\rm aw}(\bar u_{t-1}-u_{t-1})$: each step removes the fraction $k_{\rm aw}$ of the previous step's excess over the actuator limit, independently of the gains and for either sign of $K_i$. Correcting the pre-gain integral $I_t$ instead would scale the correction by $K_i$; for $K_i<0$ it becomes positive feedback, and under a constant saturating error the excess grows geometrically by the factor $1+|K_i|k_{\rm aw}$ per step. Here $k_{\rm aw}\in[0,1]$ is a per-step fraction, $k_{\rm aw}=\Delta/\tau_{\rm t}$ for a tracking time constant $\tau_{\rm t}$. With $\lambda_{\rm f}=k_{\rm aw}=0$, $v_{-1}=K_iI_{-1}$, and no residual, $v_t=K_iI_t$ recovers \eqref{eq:pid}. For $K_i=0$ the tracking state is not driven back to zero after saturation and can retain an offset; this property of back-calculation applies equally to all practical controllers. These updates remain an exact recurrent computational graph, but the three-gain quadratic convergence theorem does not apply automatically. \cref{fig:antiwindup} compares the two anti-windup laws.

\Needspace{12\baselineskip}
\section{Experimental details}\label{app:experiments}
This appendix gives the protocols of E1--E3b and E5 and the parts of the E4 protocol not stated in \cref{sec:experiments}.

\subsection{E1: exact equivalence and the LSTM distinction}
The test set had $10^4$ independent error sequences with lengths in $\{32,128,512,2048\}$, sampling periods log-uniform on $[10^{-3},1]$, gains uniform on $[-10,10]^3$, and initial memories uniform on $[-1,1]$. The errors were Gaussian white noise, PRBS ($\pm1$, switching probability $0.05$), ramps, and three-tone multisines in equal proportions, each multiplied by an error scale $s\in\{10^{-3},1,10^3\}$ so that both the absolute and the relative regime of the metric are exercised. The positional PID was compared with the two-state readout, the three-state realization \eqref{eq:threestate}, the incremental form with consistent initialization, and the hard-gated variant.

The metric was $\max_t|u_t-u'_t|/(1+\max_t|u_t|)$ for the outputs $u$ and $u'$ of two realizations in float64. Every case had to be at most $10^{-10}$; a pilot of 400 length-2048 sequences gave at most $3.9\times10^{-12}$, a margin of at least 25. Float32 errors were recorded without a threshold, and a descriptive stress run with sequences of length 65,536 measured how the rounding of the incremental form grows with the horizon. E1 is an implementation check, not an empirical contribution.

For the standard LSTM, the bound of \cref{prop:lstm} was checked on forty networks with random weights under constant error, together with the time at which a PID with $K_i\neq0$ exceeds it. A learnability probe then trained LSTMs with 8 and 32 units to output the running integral of an error carrying a random constant offset, which is the sustained-error regime for which integral action is designed. The readout was linear, from $h_t$ or from the cell state. Training horizons were $\{32,\dots,256\}$ and evaluation horizons $\{32,128,512,2048\}$, which separates interpolation from extrapolation, and the learning rate was selected on validation data so that a negative result cannot be attributed to an untuned optimizer. The probe is the only supervised component of the study, and its target is the error integral, not a PID action; E2--E5 use no gain or action labels.

\subsection{E2: recovery from physical cost only}
E2 used the plant of \cref{sec:example} with $\alpha=0$ and $H$, $R$ as specified, fresh resets $z\sim\mathcal N(0,\Id_3)$, and scalar process noise $w\sim\mathcal N(0,\sigma^2)$ entering as $\xi=b\,w$, with $\sigma\in\{0,0.1,0.5\}$. Gradient descent was run both in the gain coordinates $k$ and in the readout coordinates $\theta=\Gamma_1k$ of the recurrent implementation; an independent normal-equation solution was used only by the evaluation code, and no action labels or gain-regression losses were used. The sample size ranged over $N\in\{64,256,1024,4096,16384\}$, with $1000$ independent replications per cell. Each parameterization used $\eta=1/\lambda_{\max}$ of its own empirical Hessian, because a step valid in $k$ can diverge in $\theta$ (\cref{cor:readout}). Iterations stopped at 5,000 or at a relative gradient norm below $10^{-10}$. Both parameterizations started from the same controller, with gains uniform on $[-1,1]^3$, and a variant with gains uniform on $[-10,10]^3$ checked global convergence. A separate population run replaced the sample averages by exact expectations, computed by three-node Gauss--Hermite quadrature, which is exact for this quadratic cost; with $\Delta=1$ and $\Sigma=\Id_3$ it checks the rate $1-1/\kappa(\mathcal H_\theta)$ of \cref{cor:readout} in readout coordinates, where $\kappa(\mathcal H_\theta)=7+4\sqrt3$.

The pass rules required gradient descent to agree with the normal-equation solution to $10^{-8}$ and, for $\sigma=0$, $\|\widehat k-k_\star\|\leq10^{-8}$ for every $N$, because $u_\star$ lies in the PID span (\cref{app:finite}). For $\sigma>0$, the mean over replications of $N\|\widehat k-k_\star\|^2$ had to lie within $10\%$ of its exact value $(b^\T Hb/M)^2\sigma^2\operatorname{tr}(\Sigma^{-1})\,N/(N-4)$, which uses $\E[\widehat\Sigma^{-1}]=N\Sigma^{-1}/(N-4)$ for Gaussian resets, and the observed contraction had to match each parameterization's exact empirical rate. Two controls complete the design. Rank-deficient resets with $I=0$ served as a negative control, in which $K_i$ was required not to move from its initial value in gain coordinates, and a correlated design with $\Sigma_{ij}=0.8^{|i-j|}$ tested conditioning and had to return $k_\star$ for $\sigma=0$.

\subsection{E3: analytic nonlinear separation}
With $\alpha=1/5$, four classes were trained on the same physical cost: PID, PID plus bias, PID plus $e^2$, and PID plus bias and $e^2-1$. They are ordered by inclusion: PID lies in both PID plus bias and PID plus $e^2$, and both lie in PID plus bias and $e^2-1$. The last two contain $u_\star$. E3 used the $N$ grid and 1000 replications of E2, without process noise. Costs were computed from the plant transition, not from the completed-square identity. Analytic controllers were evaluated on one independent test set of $10^6$ Gaussian resets, paired across methods, and all controllers were also evaluated exactly by tensor Gauss--Hermite quadrature with three nodes per coordinate, which integrates the polynomial one-step cost without error. The exact evaluation was used for learned controllers, because the Monte Carlo error of the test set is common to all replications and exceeds the excess cost being measured at large $N$.

The analytic paired gaps on the test set had to agree with \eqref{eq:numericgap} and \eqref{eq:biasgap} within five Monte Carlo standard errors, the quadrature gaps exactly, and the realizable classes had to return $(k_\star,\beta_\star)$ and $(k_\star,\beta_\star,\beta_\star)$ to $10^{-8}$ in every replication (\cref{app:finite} with $\widehat\epsilon=0$). For the PID class, $\widehat k-k_\star=\beta_\star\widehat\Sigma^{-1}N^{-1}\sum_iz_ie_i^2$, and its excess population cost has leading-order mean $c_{\rm ex}/N$ with $c_{\rm ex}=M\beta_\star^2(\E e^6+\E I^2e^4+\E d^2e^4)=21M\beta_\star^2$. For PID plus bias, $c_{\rm ex}=M\beta_\star^2\,\E[(e^2-1)^2\|(z,1)\|^2]=16M\beta_\star^2$. The mean over replications had to lie within $10\%$ of $c_{\rm ex}/N$ for $N\geq1024$, so the learned gap converges to the analytic gap from above. In empirically residualized coordinates, the PID block equals $k_\star+\widehat C^\T\beta_\star$, where $\widehat C$ is the training-sample regression coefficient of the nonlinear feature on the PID features.

Several variants probe the edges of the theory. For $\alpha\in\{0,0.05,0.1,0.2,0.4\}$ the quadrature gaps are $108\alpha^2/55$ and $72\alpha^2/55$. An irrelevant feature independent of $z$ has zero population improvement. Resets $z\sim\mathcal N((1/2,0,0)^\T,\Id_3)$ test \cref{cor:raw}: the learned PID class approaches $K_p=197/550$, while raw PID plus $e^2$ returns $K_p=1/2$ exactly. Bounded resets, uniform on $[-\sqrt3,\sqrt3]^3$ with the same variance, remain within the theory, with exact gaps $324/6875$ and $144/6875$ because $\E[e^4]=9/5$. Actuator saturation lies outside the theorem and was reported without an analytic target.

\subsection{E3b: closed-loop bridge}
The plant of \cref{sec:example} was rolled out in closed loop for 200 steps from $10^5$ independent Gaussian resets. The controllers were the analytic $k_\star$, PID plus bias, and $u_\star$, the discrete-time LQR reference, and the controllers learned in 30 replications of E3 at $N=4096$. A rollout diverged if $\|z\|_\infty>10^3$ or a value became non-finite. The reported quantities are divergence fractions with Wilson $95\%$ intervals and the mean per-step cost of surviving rollouts. \Cref{prop:closedloop} predicts local exponential stability under $k_\star$ and $u_\star$; global behavior is empirical, and E3b has no pass rule. It does not test superiority.

\subsection{E4: plants, references, failure scoring, and tuning}\label{app:e4details}
\begin{table}[tbp]\centering\small
\caption{E4 benchmark families and their training parameter ranges.}\label{tab:families}
\begin{tabularx}{\textwidth}{@{}l>{\raggedright\arraybackslash}X@{}}
\toprule
Family & Dynamics and training parameter ranges\\\midrule
Linear oscillator & $\ddot y=-2\zeta\omega\dot y-\omega^2y+b_{\rm p}(\bar u+\delta_u)+w$; $\zeta\in[0.2,0.6]$, $\omega\in[0.8,1.2]$, $b_{\rm p}\in[0.8,1.2]$.\\
Duffing oscillator & Same, with additional drift $-\gamma y^3$; $\gamma\in[0.1,0.5]$ and the same remaining ranges.\\
Hammerstein & $x_{t+1}=a_{\rm p}x_t+b_{\rm p}\tanh(c_{\rm p}(\bar u_t+\delta_u))+w_t$, $y_t=x_t$; $a_{\rm p}\in[0.85,0.95]$, $b_{\rm p}\in[0.2,0.3]$, $c_{\rm p}\in[0.8,1.2]$.\\\bottomrule
\end{tabularx}
\end{table}

The two oscillators are continuous-time plants, simulated by RK4 with four substeps and the control held fixed over $\Delta=0.02$; the Hammerstein plant is discrete, with $\Delta=1$ for the PID features. Episodes last $T=1000$ steps, with scales $s_e=1$ and $s_u=2$ and actuator limit $u_{\max}=2$. From step 500 onward a constant actuator offset $\delta_u$ acts on the plant, drawn uniformly from $[-0.2,0.2]$ for the oscillators and from $[-0.02,0.02]$ for Hammerstein. The oscillators receive white acceleration noise held over each step, with a standard deviation drawn per episode from $U[0,0.05]$; Hammerstein receives $w_t\sim\mathcal N(0,0.005^2)$. Observation noise is independent, with a standard deviation drawn per episode from $U[0,0.02]$. Episodes start with $y$ or $x$ uniform in $[-0.2,0.2]$, oscillator velocity uniform in $[-0.1,0.1]$, and $v_{-1}=0$, $D_{-1}=0$, $e_{-1}=e_0$.

The parameter panels use scrambled Sobol sequences with seeds derived from the family and split. The shift panel draws each parameter from the two outer intervals adjacent to its training range, each 20\% of that range wide, with reference amplitudes up to $1.2$ and an observation-noise standard deviation of $0.04$. Before any controller was evaluated, a fixed steady-state rule computed from parameters and references labeled the (plant, reference) pairs that $u_{\max}$ cannot hold.

\begin{figure}[t]
\centering
\includegraphics[width=0.62\textwidth]{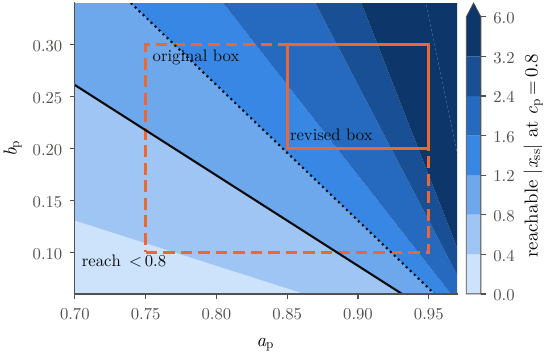}
\caption{Reachable steady-state output of the Hammerstein family under $u_{\max}$ at the least favorable input gain $c_{\rm p}$, including the disturbance, with the training box of $a_{\rm p}$ and $b_{\rm p}$ before (dashed box) and after (solid box) it was narrowed during protocol design so that every training reference is reachable. Solid contour: the training reference amplitude; dotted contour: the largest shift-test amplitude. Analytic computation; no experimental data.}
\label{fig:feasibility}
\end{figure}

Training references were an equal mixture of four types: piecewise-constant amplitudes in $[-0.8,0.8]$ changing every 200 steps, linear ramps with slopes $\pm0.002$ per step clipped at $\pm0.8$, PRBS of amplitude $\pm0.8$ switching with probability $0.02$, and multisines with frequencies $\{0.001,0.002,0.004\}$ cycles per step, random phases, and component amplitudes $0.8/3$. Parameter ranges and references were chosen so that holding a reference is feasible under $u_{\max}$. On the training ranges, the worst-case steady-state requirement is $1.96$ for constant references, on the Duffing oscillator, and $1.74$ for multisines on the linear oscillator. For multisines on the Duffing oscillator, a conservative bound reaches $2.06$ only at the extreme corner of the ranges. \Cref{fig:feasibility} shows the constraint for the Hammerstein family.

A rollout was aborted if any physical state exceeded magnitude 20 or a value became non-finite. An aborted episode scores its accumulated cost plus the remaining steps at the per-step cost of the abort threshold, $c_{\rm cap}=(20/s_e)^2+10^{-2}(u_{\max}/s_u)^2+10^{-2}(2u_{\max}/s_u)^2=400.05$, with its terminal term also at $c_{\rm cap}$, so episode scores are bounded by $800.10$. The failure rate is a co-primary endpoint, because a lower cost among surviving runs cannot by itself establish superiority. A score of $10^6$ per failure was kept as a sensitivity analysis; no episode failed (\cref{sec:results}), so neither choice affects any result.

The learned controllers' configuration grid used the learning rates $\{10^{-3},3\cdot10^{-3},10^{-2}\}$. The warm-started residual trained all parameters from the start, and the compute spent on its envelope is reported separately (\cref{sec:results}); validation computations are counted separately for every method.

The envelope searched the box $K_p,K_i,K_d\in[-10,10]$, $\lambda_{\rm f}\in[0,0.99]$, $k_{\rm aw}\in[0,1]$, extended by $\beta_0\in[-2,2]$ for PID with bias. It evaluated 4,096 scrambled Sobol candidates, refined the best 32 by bounded gradient optimization, and ran differential evolution with 10 seeds and up to 20,000 candidate evaluations each. When a selected value lay within 1\% of a box face, its interval was doubled once and the search repeated, and a value still on a face labeled the envelope boundary-limited.

\subsection{E4: precision pilot, secondary metrics, and ablation}
Before the freeze, a pilot on validation draws only, with three seeds per method, estimated the seed, plant, and episode variance components and projected the half-width of the bound on $\mathcal R_f$. The variation between episodes within a plant dominated and the seed components were negligible, so the frozen design kept 20 seeds and set 64 test episodes per plant. A rerun of the pilot after the freeze led amendment 1 to raise this to 128, which brings the projected half-width below $0.02$ in every family.

The secondary metrics were the integrated absolute error (IAE), its time-weighted version (ITAE), control energy, total variation, saturation fraction, and training cost, and learned gains, residual magnitude, and gradient norms were logged as diagnostics. Empirical failures in nonlinear rollouts are reported without being presented as stability proofs. Inference latency, overshoot and settling time on step references, and closed-loop eigenvalues for linear unsaturated cases were specified by the frozen protocol but not computed.

The frozen protocol also required one-at-a-time ablations of the nonlinear residual, the integral and derivative features, the PID warm start, and the anti-windup, which amendment 3 dropped. The two residual arms differ only in how their PID part starts, at generic gains with a PID-only phase or at the envelope, so their comparison still isolates the envelope warm start. E4 was to be interpreted only if E1--E3 passed, which they did.

\subsection{E4: the exploratory attribution arm}\label{app:attribution}
The arm of \cref{sec:attribution} trains the practical PID alone, its five parameters started at the envelope of each family (the same file, by hash, as the warm-started residual), for 20,000 episodes in batches of 16 with Adam and gradient clipping at norm 10, with 20 seeds and no PID-only warm-up. Training episodes depend only on the family, the seed, and the update, so each seed sees exactly the training episodes of the same seed of the warm-started residual. The learning rate was selected on validation from the residual's three values: $10^{-2}$, the largest, for the linear family, and $10^{-3}$ for the Duffing and Hammerstein families. The Hammerstein run used the GPU, as the other E4 runs did, and the oscillator run the CPU, for memory reasons; in float64 the two devices agree to rounding. An interrupted first attempt of the oscillator run was discarded and the run restarted from the beginning (the repository logs every attempt). A test checks that the arm follows, to float64 precision, the PID parameters of a warm-started residual whose residual is switched off.

The evaluation regenerated the test and shift episodes of the unsealing run from their seeds, checked them by plant, episode, and seed, and scored only the new arm. The analysis uses the hierarchical bootstrap of the primary endpoint with 10,000 replicates, resampling plants, episodes within plants, and seeds, with the seeds paired between the two learned arms, and reports one-sided bounds at level $1-0.05/3$. The rule for reading it was recorded before the arm was trained: an upper bound of $\mathcal R$ for the residual against the PID tuned alone below zero means that the recurrent part lowers the cost beyond PID-only tuning in that family; a lower bound above $-0.01$ means that any contribution is below one percent; otherwise the data do not separate the two. The first two readings can hold together, as in the Hammerstein family on test, where the recurrent part lowers the cost by less than one percent.

\subsection{E5: classical plants}\label{app:e5}
E5 uses eight classical PID benchmarks, all in deviation coordinates around their operating point (\cref{tab:e5plants}). The continuous-time plants hold the applied action over each sampling period and are integrated by RK4 with four substeps; the thermal chamber is discretized exactly. The series RLC circuit $L\ddot q+R\dot q+q/C=v$ is the mass-spring equation with $\omega=1/\sqrt{LC}$, $\zeta=(R/2)\sqrt{C/L}$ and $b=1/L$. The double and double inverted pendulums were left out: they are underactuated, and a single-loop PID cannot stabilize them. Each episode lasts 600 steps, with a reference step at the start and a load step at mid-episode.

\begin{table}[tbp]\centering\small
\caption{E5 plants, with sampling period $\Delta$, actuator limit, and the sizes of the reference and load steps.}\label{tab:e5plants}
\begin{tabularx}{\textwidth}{@{}l>{\raggedright\arraybackslash}Xcccc@{}}
\toprule
plant & model & $\Delta$ & limit & reference & load\\\midrule
mass-spring / RLC & $\ddot y=-2\zeta\omega\dot y-\omega^2y+bu$, $\omega=b=1$, $\zeta\in\{0.2,1,2\}$ & 0.05 & 2 & 1 & 0.2\\\rowgap
pendulum & $\ddot\vartheta=-\sin\vartheta-0.1\dot\vartheta+u$ & 0.05 & 2 & 1 & 0.2\\\rowgap
inverted pendulum & $\ddot\vartheta=\sin\vartheta-0.1\dot\vartheta+u$, torque at the pivot & 0.05 & 2 & 0.2 & 0.1\\\rowgap
coupled tanks & $\dot h_1=1+u-\sqrt{h_1}$, $\dot h_2=\sqrt{h_1}-\sqrt{h_2}$, output $h_2$ & 0.1 & 1 & 0.3 & $-0.1$\\\rowgap
thermal chamber & $\dot T=-T+u(t-L)$, $L\in\{0.1,0.3\}$ & 0.02 & 1 & 0.5 & $-0.2$\\\bottomrule
\end{tabularx}
\end{table}

E5a compares, on every plant, the textbook PID laws, written independently of the library code, with the networks that realize them, each in its own closed loop with the same reference, load, and noise. The ideal PID \eqref{eq:pid} is compared with the two-state network of \cref{th:real} and with the incremental form, the practical PID of \cref{sec:long} with the E4 controller and its linear recurrent cell, and the E4 residual controller with its output layer at zero with its PID. The gain sets are a classical tuning and 20 random draws around it. The classical tuning places the poles of the second-order linearization at $(s+p)^3$, with $p=2$ for the mechanical plants and $p=1$ for the tanks, and uses the SIMC rule~\citep{skogestad2003} with $\tau_c=L$ for the thermal chamber; its practical parameters are $\lambda_{\rm f}=0.5$ and $k_{\rm aw}=0.1$. The random draws multiply each gain by a log-uniform factor in $[1/3,3]$ and draw $\lambda_{\rm f}$ and $k_{\rm aw}$ uniformly, and a draw is kept only if both textbook laws keep every scenario bounded. The three scenarios are a loop without clipping or noise, a loop with clipping at the actuator limit, and the latter with observation and process noise. The pass rule bounds the normalized discrepancy $\max_t|a_t-b_t|/(1+\max_t|b_t|)$ of the applied and commanded actions and of the output by $10^{-10}$, as in E1.

E5b trains the network on the one-step physical cost of the six plants without dead time; with a dead time the action does not reach the next measured state, and the one-step cost does not identify the gains. The one-step model acts on $z=(x_1,x_2,I,e_{t-1})$ with reference zero and is affine in the action, as (A1) requires: semi-implicit Euler for the mechanical plants, so that the action reaches the output within one step, and explicit Euler for the tanks. The cost weights are $H=\operatorname{diag}(1,0.1,0.1,0)$, or $\operatorname{diag}(0.1,1,0.1,0)$ for the tanks, and $R=10^{-4}$, on 4,096 Gaussian resets with $e_{t-1}$ set so that $D_t=-\dot y$. Gradient descent starts from zero gains, in gain and in readout coordinates, with step $1/\lambda_{\max}$ of its Hessian, and stops at a relative gradient norm of $10^{-14}$; the normal equations are computed on the evaluation side only. The pass rules require the learned gains to match the normal equations to $10^{-8}$ relative to $1+\|k_\star\|$ in both coordinates, the contraction to stay within the predicted rate, and the closed loop of the trained network to match that of the PID with the normal-equation gains to $10^{-6}$ wherever that loop is bounded. A pilot fixed $R$ and the stopping tolerance before the reported runs. The thresholds for the realizations and for the learned gains are those of E1 and E2, and the closed-loop threshold, $10^{-6}$, was fixed with them.

\subsection{Amendments after the freeze}\label{sec:amend}
All amendments were dated 21 October 2026 and recorded with their motivation in the repository's protocol document, which labels them A1--A4. Amendment 1 gave each test and shift plant 128 episodes instead of the 64 fixed at the freeze, because the variation between episodes within a plant dominated the pilot. Amendment 2 restricted the gradient-trained PID in readout coordinates (the recurrent PID) to the configuration selected for the same PID trained in gain coordinates (the direct PID), and amendment 3 later dropped both. Amendment 3 reduced E4, for compute, to the primary endpoint and the five controllers of \cref{sec:experiments}, dropping the other learned controllers, the ablations, and the oracle PID. Amendment 4 gave the decision rule explicit outcome classes; it was committed while the single test evaluation was running and before any of its output was read.

One correction was also made. The first envelope runs selected only among the finalists of the doubled search box, discarding first-box finalists that validated better. Selection was corrected to range over both boxes, and every envelope and the warm-started residual were rerun before the test panels were generated. The search itself was deterministic: the candidates were bitwise identical, and only the selected point changed.

\section{Additional results}\label{app:results}
This appendix gives the full results of E1--E3b and E5, summarized in \cref{sec:results}, and the learning curves and secondary metrics of E4.

\subsection{E1: exact equivalence and the LSTM distinction} 
\begin{figure}[tbp]
\centering
\includegraphics[width=\textwidth]{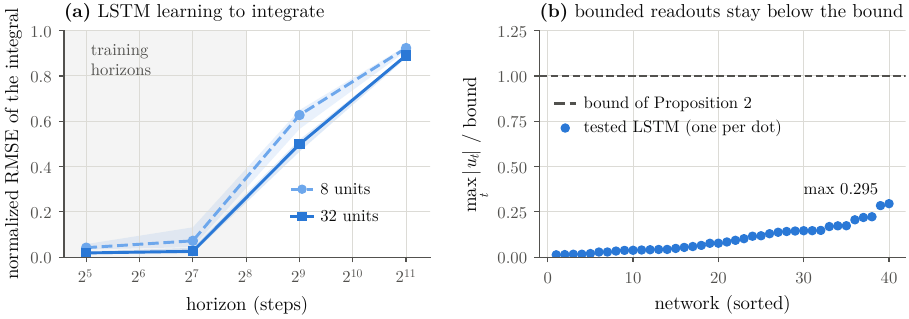}
\caption{E1: (a) normalized error of LSTMs trained to output the running integral, by evaluation horizon (median over 5 seeds with their range shaded; the gray band marks the training horizons); (b) the largest output of each tested LSTM relative to the bound of \cref{prop:lstm}.}
\label{fig:e1}
\end{figure}

All \EoneCases{} cases meet the $10^{-10}$ tolerance, with a maximum normalized discrepancy of $\EoneMax$ (\cref{tab:e1}). In float32 the incremental form drifts to $3.7\times10^{-4}$ and the other realizations to $7.1\times10^{-6}$. In the float64 stress run, the discrepancy of the incremental form grows from $4.0\times10^{-14}$ at length 2048 to $4.7\times10^{-13}$ at length 65,536. The LSTM bound of \cref{prop:lstm} held for all 40 networks tested, and a PID with $K_i=1$ exceeded it after 56 to 897 steps. The learnability probe fits the running integral within the training horizons, with a median normalized RMSE of 0.02--0.07 at horizons up to 128. It fails beyond them, with 0.89--0.92 at horizon 2048 (\cref{fig:e1}).

\begin{table}[tbp]\centering\small
\caption{E1: normalized action discrepancy between realizations, over all cases (float64, and float32 maximum). Green: within the $10^{-10}$ tolerance of the pass rule.}\label{tab:e1}
\begin{tabular}{lrrrrr}
\toprule
realization & median & 90\% & 99\% & max & max (float32) \\
\midrule
two-state readout & $2.1\times10^{-16}$ & $1.3\times10^{-15}$ & $1.1\times10^{-14}$ & \cpass $2.3\times10^{-14}$ & $7.1\times10^{-6}$ \\
three-state & $2.1\times10^{-16}$ & $1.3\times10^{-15}$ & $1.1\times10^{-14}$ & \cpass $2.3\times10^{-14}$ & $7.1\times10^{-6}$ \\
hard-gated cell & $2.1\times10^{-16}$ & $1.3\times10^{-15}$ & $1.1\times10^{-14}$ & \cpass $2.3\times10^{-14}$ & $7.1\times10^{-6}$ \\
incremental & $1.0\times10^{-15}$ & $1.2\times10^{-14}$ & $6.2\times10^{-14}$ & \cpass $6.2\times10^{-13}$ & $3.7\times10^{-4}$ \\
\bottomrule
\end{tabular}

\end{table}

\subsection{E2: recovery from physical cost only} 
\begin{figure}[tbp]
\centering
\includegraphics[width=\textwidth]{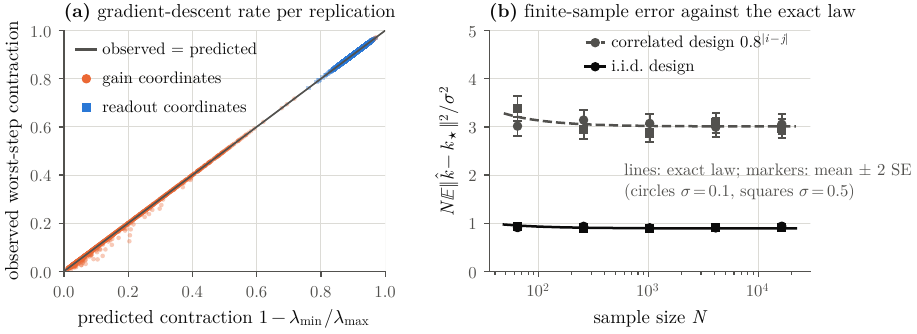}
\caption{E2: (a) observed worst-step contraction of gradient descent (GD) against its predicted value, per replication, in gain and readout coordinates; (b) mean scaled gain error $N\|\widehat k-k_\star\|^2/\sigma^2$ against its exact finite-sample law, with $\pm2$ standard errors, for the i.i.d.\ and correlated designs.}
\label{fig:e2}
\end{figure}

Gradient descent agrees with the normal equations to $\EtwoNEMaxCorr$, and noiseless recovery is exact to $10^{-8}$ for every $N$, in both designs. Over \EtwoReps{} replications per cell, the noise law of \cref{app:finite} holds within the 10\% rule. The largest deviation from the exact value is $|z|=\EtwoMaxAbsZ$ standard errors (\cref{fig:e2,tab:e2}). The readout coordinates contract at the predicted rate $\EtwoReadoutRate$. In the rank-deficient negative control, $K_i$ stays at its initial value in gain coordinates.

\begin{table}[htbp]\centering\small
\caption{E2: mean over 1000 replications of the scaled gain error $N\|\widehat k-k_\star\|^2$ against its exact finite-sample value, for the i.i.d.\ and correlated designs. Green ratios lie within the $\pm10\%$ pass rule.}\label{tab:e2}
\begin{tabular}{lrrrrrrr}
\toprule
design & $N$ & $\sigma$ & mean $N\|\hat k-k_\star\|^2$ & SE & exact law & ratio & $z$ \\
\midrule
i.i.d. & 64 & 0.1 & 0.00909 & 0.00024 & 0.00952 & \cpass 0.954 & $-$1.82 \\
 &  & 0.5 & 0.23308 & 0.00639 & 0.23802 & \cpass 0.979 & $-$0.77 \\
\addlinespace[2pt]
 & 256 & 0.1 & 0.00936 & 0.00027 & 0.00907 & \cpass 1.032 & $+$1.10 \\
 &  & 0.5 & 0.22229 & 0.00602 & 0.22668 & \cpass 0.981 & $-$0.73 \\
\addlinespace[2pt]
 & 1024 & 0.1 & 0.00906 & 0.00023 & 0.00896 & \cpass 1.011 & $+$0.42 \\
 &  & 0.5 & 0.22314 & 0.00590 & 0.22402 & \cpass 0.996 & $-$0.15 \\
\addlinespace[2pt]
 & 4096 & 0.1 & 0.00903 & 0.00023 & 0.00893 & \cpass 1.010 & $+$0.41 \\
 &  & 0.5 & 0.22725 & 0.00586 & 0.22336 & \cpass 1.017 & $+$0.66 \\
\addlinespace[2pt]
 & 16384 & 0.1 & 0.00938 & 0.00024 & 0.00893 & \cpass 1.051 & $+$1.92 \\
 &  & 0.5 & 0.22814 & 0.00589 & 0.22319 & \cpass 1.022 & $+$0.84 \\
\midrule
$0.8^{|i-j|}$ & 64 & 0.1 & 0.03017 & 0.00099 & 0.03209 & \cpass 0.940 & $-$1.93 \\
 &  & 0.5 & 0.84505 & 0.03183 & 0.80220 & \cpass 1.053 & $+$1.35 \\
\addlinespace[2pt]
 & 256 & 0.1 & 0.03142 & 0.00105 & 0.03056 & \cpass 1.028 & $+$0.82 \\
 &  & 0.5 & 0.73588 & 0.02422 & 0.76400 & \cpass 0.963 & $-$1.16 \\
\addlinespace[2pt]
 & 1024 & 0.1 & 0.03073 & 0.00101 & 0.03020 & \cpass 1.018 & $+$0.53 \\
 &  & 0.5 & 0.71816 & 0.02334 & 0.75502 & \cpass 0.951 & $-$1.58 \\
\addlinespace[2pt]
 & 4096 & 0.1 & 0.02994 & 0.00097 & 0.03011 & \cpass 0.994 & $-$0.18 \\
 &  & 0.5 & 0.77021 & 0.02544 & 0.75280 & \cpass 1.023 & $+$0.68 \\
\addlinespace[2pt]
 & 16384 & 0.1 & 0.03055 & 0.00104 & 0.03009 & \cpass 1.015 & $+$0.45 \\
 &  & 0.5 & 0.74135 & 0.02406 & 0.75225 & \cpass 0.986 & $-$0.45 \\
\bottomrule
\end{tabular}

\end{table}

\subsection{E3 and E3b: analytic separation and closed-loop bridge} 
\begin{figure}[tbp]
\centering
\includegraphics[width=\textwidth]{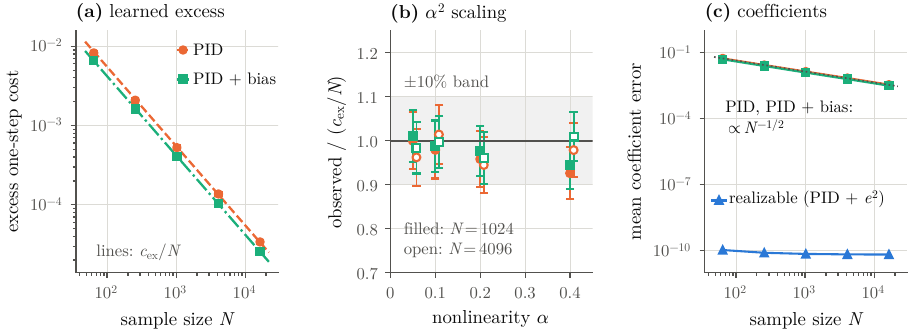}
\caption{E3: (a) learned excess one-step cost against $c_{\rm ex}/N$; (b) observed over predicted excess across the nonlinearity $\alpha$, with the $\pm10\%$ band of the pass rule; (c) mean coefficient error, which decays as $N^{-1/2}$ for the PID classes and stays at round-off for the realizable class.}
\label{fig:e3}
\end{figure}

The exact quadrature gaps equal $108\alpha^2/55$ and $72\alpha^2/55$ to machine precision for every $\alpha$, and the paired Monte Carlo gaps lie within 1.35 standard errors of them (\cref{tab:e3}). At $\alpha=1/5$ the gap is \GapRelPid{} of $\mathcal L_\PID^\star=\LpidStar$ against PID and \GapRelPidBias{} against PID with bias. The realizable classes recover $(k_\star,\beta_\star)$ to $\EthreeRealizableMax$ in all \EthreeReps{} replications. For $N\geq1024$ the learned excess cost lies within 2\% of $c_{\rm ex}/N$ at $\alpha=1/5$ and within the 10\% rule for every $\alpha$, and the log--log slopes are $-0.99$ (\cref{fig:e3}). In closed loop (\cref{tab:e3b}), the optimal PID diverged in \EthreebDivPid{} of the Gaussian resets. PID with bias diverged in \EthreebDivPidBias{}, with a 16-fold higher cost among survivors. The one-step optimum $u_\star$ diverged in none of $10^5$ resets.

\begin{table}[tbp]\centering\small
\caption{E3: one-step cost gaps (closed form, exact quadrature, independent NumPy path, paired Monte Carlo). Green entries meet their pass rule: exact agreement, or Monte Carlo within five standard errors.}\label{tab:e3}
\begin{tabular}{llrrrrr}
\toprule
$\alpha$ & baseline & closed form & quadrature & independent path & Monte Carlo & SE \\
\midrule
0.05 & PID & 0.004909 & \cpass 0.004909 & \cpass 0.004909 & \cpass 0.004928 & $1.6\times10^{-5}$ \\
 & PID + bias & 0.003273 & \cpass 0.003273 & \cpass 0.003273 & \cpass 0.003289 & $1.2\times10^{-5}$ \\
\addlinespace[2pt]
0.1 & PID & 0.019636 & \cpass 0.019636 & \cpass 0.019636 & \cpass 0.019714 & $6.5\times10^{-5}$ \\
 & PID + bias & 0.013091 & \cpass 0.013091 & \cpass 0.013091 & \cpass 0.013158 & $5.0\times10^{-5}$ \\
\addlinespace[2pt]
0.2 & PID & 0.078545 & \cpass 0.078545 & \cpass 0.078545 & \cpass 0.078856 & $2.6\times10^{-4}$ \\
 & PID + bias & 0.052364 & \cpass 0.052364 & \cpass 0.052364 & \cpass 0.052631 & $2.0\times10^{-4}$ \\
\addlinespace[2pt]
0.4 & PID & 0.314182 & \cpass 0.314182 & \cpass 0.314182 & \cpass 0.315422 & $1.0\times10^{-3}$ \\
 & PID + bias & 0.209455 & \cpass 0.209455 & \cpass 0.209455 & \cpass 0.210524 & $7.9\times10^{-4}$ \\
\bottomrule
\end{tabular}

\end{table}

\begin{table}[tbp]\centering\small
\caption{E3b: closed-loop divergence over 200 steps from $10^5$ Gaussian resets, and mean cost of surviving rollouts; the last column gives the mean [minimum, maximum] over the 30 learned fits. There is no pass rule; green marks the lowest value in each column.}\label{tab:e3b}
\begin{tabular}{lrrrr}
\toprule
controller & diverged (\%) & Wilson 95\% & cost (survivors) & learned, 30 fits (\%) \\
\midrule
PID $k_\star$ & 0.524 & [0.481, 0.571] & 0.0341 & 0.53 [0.42, 0.63] \\
PID + bias & 0.300 & [0.268, 0.336] & 0.5465 & 0.30 [0.24, 0.36] \\
$u_\star$ (PID + $e^2$) & \cpass 0.000 & [0.000, 0.004] & 0.0181 & \cpass 0.00 [0.00, 0.00] \\
LQR (reference) & 0.033 & [0.023, 0.046] & \cpass 0.0153 & -- \\
\bottomrule
\end{tabular}

\end{table}

\subsection{E5: classical plants} 
\begin{figure}[tbp]
\centering
\includegraphics[width=\textwidth]{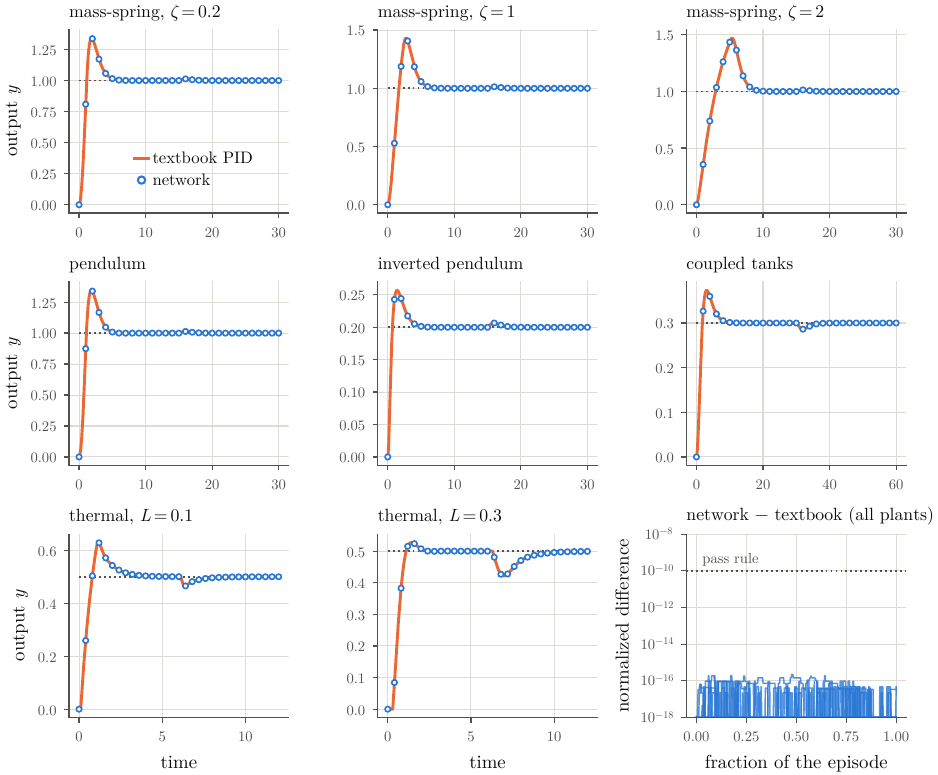}
\caption{E5: closed-loop responses of the textbook PID (line) and of the network of \cref{th:real} (markers) on the eight classical plants, as in \cref{fig:e5}; the last panel repeats the normalized output difference for every plant.}
\label{fig:e5all}
\end{figure}

\begin{figure}[tbp]
\centering
\includegraphics[width=\textwidth]{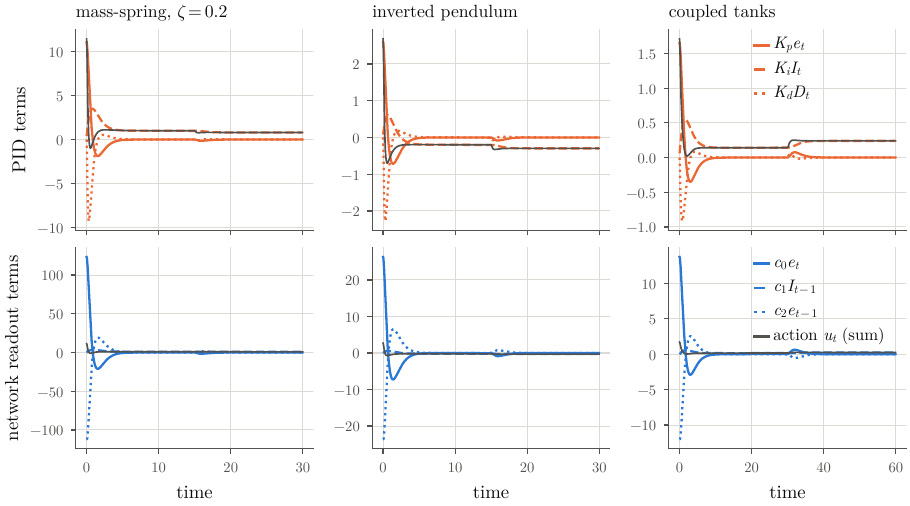}
\caption{E5: the same action decomposed by the PID (top: proportional, integral and derivative terms) and by the network (bottom: its readout terms), with the classical tuning and no clipping. The sums agree; the network's terms are about an order of magnitude larger and cancel.}
\label{fig:e5terms}
\end{figure}

\begin{figure}[tbp]
\centering
\includegraphics[width=\textwidth]{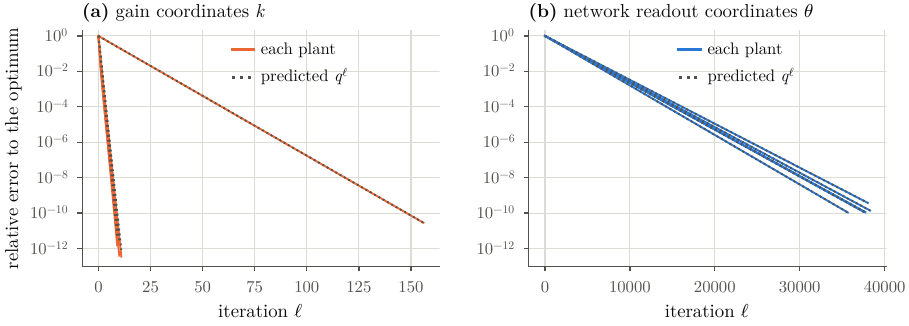}
\caption{E5b: error of gradient descent on the one-step cost, relative to its initial value, in gain coordinates (a) and in the network's readout coordinates (b), for the six plants, with the predicted contraction $q^\ell$ dotted. In gain coordinates all six converge within \EfiveItersGain{} iterations, so their curves overlap.}
\label{fig:e5conv}
\end{figure}

\begin{figure}[tbp]
\centering
\includegraphics[width=\textwidth]{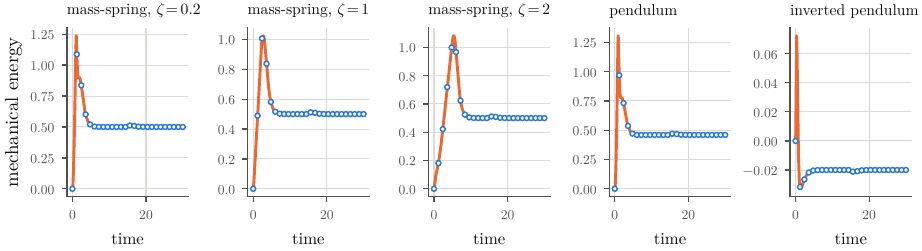}
\caption{E5: mechanical energy along the closed loops of the textbook PID (line) and the network (markers) for the five mechanical plants, with the classical tuning and actuator saturation.}
\label{fig:e5energy}
\end{figure}

\begin{figure}[tbp]
\centering
\includegraphics[width=\textwidth]{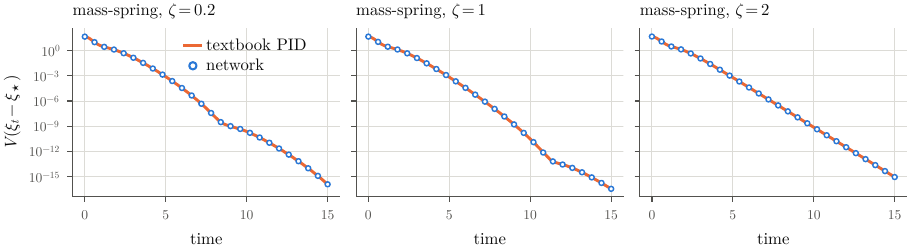}
\caption{E5: a quadratic Lyapunov function of the closed loop along the textbook PID (line) and the network (markers), on the three mass-spring plants with the classical tuning, without clipping, from the reference step to the load step.}
\label{fig:e5lyap}
\end{figure}

Every pass rule was met (\cref{tab:e5a,tab:e5b}), and \cref{fig:e5all} shows the closed loops on all eight plants. The practical PID of E4 agrees exactly with its textbook implementation, because both perform the same operations in the same order, and its residual controller with the output layer at zero is bitwise identical to it. \Cref{fig:e5terms} shows the two decompositions of the same action, \cref{fig:e5conv} the convergence of gradient descent in both coordinates, and \cref{fig:e5energy} the mechanical energy of both loops. \Cref{fig:e5lyap} shows the stabilization through a Lyapunov function. On the mass-spring plants without clipping, the closed loop of the ideal PID is linear in the deviation $\xi_t-\xi_\star$ of its state $\xi_t=(y_t,\dot y_t,I_{t-1},e_{t-1})$ from the equilibrium, so $V(\xi)=\xi^\T P\xi$, with $P$ from the discrete Lyapunov equation of that loop and $Q=\Id_4$, decreases at every step; along the textbook PID and along the network it agrees to $\LyapEfiveRel$ of its initial value.

\begin{table}[tbp]\centering\footnotesize\setlength{\tabcolsep}{4pt}
\caption{E5a: largest normalized discrepancy in closed loop between each network and its textbook PID, over all gain sets and scenarios.}\label{tab:e5a}
\begin{tabular}{lccccc}
\toprule
 & \multicolumn{2}{c}{ideal PID} & \multicolumn{2}{c}{practical PID} & warm start \\
\cmidrule(lr){2-3}\cmidrule(lr){4-5}
plant & two-state network & incremental form & E4 practical PID & E4 recurrent cell & (residual at zero) \\
\midrule
mass-spring, $\zeta=0.2$ & $2\times10^{-13}$ & $4\times10^{-12}$ & 0 & $5\times10^{-14}$ & identical \\
mass-spring, $\zeta=1$ & $4\times10^{-14}$ & $2\times10^{-13}$ & 0 & $4\times10^{-14}$ & identical \\
mass-spring, $\zeta=2$ & $1\times10^{-14}$ & $8\times10^{-14}$ & 0 & $5\times10^{-15}$ & identical \\
pendulum & $1\times10^{-13}$ & $6\times10^{-13}$ & 0 & $6\times10^{-14}$ & identical \\
inverted pendulum & $2\times10^{-14}$ & $4\times10^{-14}$ & 0 & $2\times10^{-14}$ & identical \\
coupled tanks & $7\times10^{-14}$ & $5\times10^{-13}$ & 0 & $7\times10^{-15}$ & identical \\
thermal, $L=0.1$ & $5\times10^{-15}$ & $3\times10^{-14}$ & 0 & $3\times10^{-15}$ & identical \\
thermal, $L=0.3$ & $2\times10^{-15}$ & $3\times10^{-14}$ & 0 & $2\times10^{-15}$ & identical \\
\bottomrule
\end{tabular}

\end{table}

\begin{table}[tbp]\centering\footnotesize\setlength{\tabcolsep}{4pt}
\caption{E5b: one-step optimal gains from the normal equations, error of the gains learned by gradient descent in gain coordinates and in the network's readout coordinates (relative to $1+\|k_\star\|$), iterations, condition number of the readout Hessian, and closed-loop discrepancy (output and applied action) between the trained network and the PID with the normal-equation gains.}\label{tab:e5b}
\begin{tabular}{lrrrcccccc}
\toprule
 & \multicolumn{3}{c}{one-step optimal gains} & \multicolumn{2}{c}{error of the learned gains} & \multicolumn{2}{c}{iterations} & & closed \\
\cmidrule(lr){2-4}\cmidrule(lr){5-6}\cmidrule(lr){7-8}
plant & $K_p$ & $K_i$ & $K_d$ & gains & network & gains & network & $\kappa(\mathcal H_\theta)$ & loop \\
\midrule
mass-spring, $\zeta=0.2$ & 6.300 & 0.0351 & 14.098 & $3\times10^{-12}$ & $1\times10^{-10}$ & 10 & 37{,}761 & 1{,}649 & $8\times10^{-11}$ \\
mass-spring, $\zeta=1$ & 6.300 & 0.0351 & 12.947 & $3\times10^{-13}$ & $1\times10^{-10}$ & 11 & 37{,}498 & 1{,}642 & $7\times10^{-11}$ \\
mass-spring, $\zeta=2$ & 6.300 & 0.0351 & 11.509 & $2\times10^{-12}$ & $1\times10^{-10}$ & 10 & 38{,}286 & 1{,}688 & $7\times10^{-11}$ \\
pendulum & 6.384 & 0.0337 & 14.314 & $4\times10^{-13}$ & $9\times10^{-11}$ & 10 & 35{,}668 & 1{,}555 & $8\times10^{-11}$ \\
inverted pendulum & 7.649 & 0.0348 & 14.313 & $1\times10^{-12}$ & $9\times10^{-11}$ & 9 & 37{,}704 & 1{,}645 & $5\times10^{-11}$ \\
coupled tanks & 8.565 & 0.0249 & 16.872 & $3\times10^{-11}$ & $4\times10^{-10}$ & 156 & 38{,}025 & 1{,}753 & $3\times10^{-10}$ \\
\bottomrule
\end{tabular}

\end{table}

\subsection{E4: relative costs, learning curves, and secondary metrics} 
\begin{figure}[tbp]
\centering
\includegraphics[width=\textwidth]{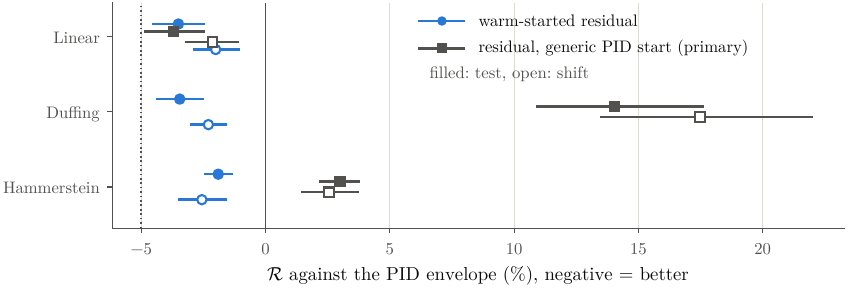}
\caption{E4: relative cost against the PID envelope, with one-sided bounds at level $1-0.05/3$, for test (filled) and shift (open) panels. The dotted line is the margin of the primary endpoint.}
\label{fig:e4}
\end{figure}

\begin{figure}[tbp]
\centering
\includegraphics[width=\textwidth]{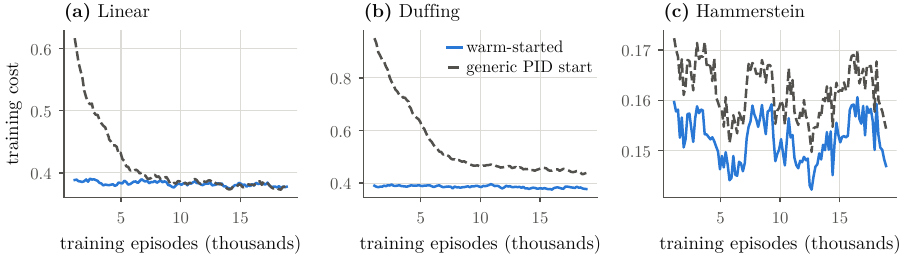}
\caption{E4: training cost of the two residual arms (median over the 20 seeds of the selected configuration, smoothed by a centered window of 15 logged updates). Both arms see the same training episodes in each update, so the fluctuations are shared.}
\label{fig:e4learning}
\end{figure}

\Cref{fig:e4} plots the relative costs of \cref{tab:e4}, \cref{fig:e4learning} shows the training cost of the two residual arms, and \cref{tab:e4sec} the secondary metrics on the test panel.

\begin{table}[tbp]\centering\small
\caption{E4: test means of the episode score and its secondary metrics (IAE and ITAE of the scaled error, control energy, total variation of the applied action, fraction of saturated steps). All are costs; green marks the lowest value per family.}\label{tab:e4sec}
\begin{tabular}{llrrrrrr}
\toprule
family & controller & score & IAE & ITAE & energy & TV & saturated (\%) \\
\midrule
Linear & Warm-started residual & 0.3769 & 0.275 & 0.134 & \cpass 0.675 & 235.36 & 8.5 \\
 & Residual, generic PID start & \cpass 0.3762 & 0.268 & 0.132 & 0.712 & 153.99 & 7.3 \\
 & PID envelope & 0.3906 & 0.271 & 0.129 & 0.884 & 186.74 & 7.3 \\
 & PID-plus-bias envelope & 0.3916 & 0.277 & 0.131 & 0.840 & \cpass 141.94 & \cpass 4.7 \\
 & Equal-budget PID & 0.4556 & \cpass 0.220 & \cpass 0.105 & 1.254 & 272.04 & 13.5 \\
\midrule
Duffing & Warm-started residual & \cpass 0.3823 & 0.273 & 0.133 & 0.731 & 231.49 & 9.2 \\
 & Residual, generic PID start & 0.4515 & 0.310 & 0.156 & \cpass 0.582 & 437.50 & \cpass 1.9 \\
 & PID envelope & 0.3960 & 0.272 & 0.130 & 0.876 & 243.88 & 7.8 \\
 & PID-plus-bias envelope & 0.3973 & 0.294 & 0.141 & 1.304 & 379.61 & 18.0 \\
 & Equal-budget PID & 0.4274 & \cpass 0.241 & \cpass 0.116 & 1.150 & \cpass 108.65 & 15.1 \\
\midrule
Hammerstein & Warm-started residual & \cpass 0.1548 & 0.074 & 0.037 & 0.169 & 133.82 & 2.2 \\
 & Residual, generic PID start & 0.1626 & \cpass 0.051 & \cpass 0.025 & \cpass 0.153 & \cpass 55.98 & \cpass 1.3 \\
 & PID envelope & 0.1578 & 0.057 & 0.028 & 0.187 & 118.93 & 2.6 \\
 & PID-plus-bias envelope & 0.1583 & 0.062 & 0.031 & 0.189 & 122.85 & 2.6 \\
 & Equal-budget PID & 0.1623 & 0.061 & 0.030 & 0.184 & 174.06 & 2.3 \\
\bottomrule
\end{tabular}

\end{table}

\section{Notation}\label{app:notation}
Transposes are written $(\cdot)^\T$ in sans-serif, distinct from the horizon $T$; $\Id_n$ is the $n\times n$ identity and $\E$ an expectation. Symbols with a subscript are distinct from the bare letter: for example $b_{\rm p}$ (a plant gain) differs from the input vector $b$, and $s_e,s_u$ (scales) differ from the state $s_t$. The plants of E5 (\cref{app:e5}) keep their conventional symbols locally: $\vartheta$ for a pendulum angle, $h_1,h_2$ for tank levels, $T$ for the chamber temperature, $L$ for its dead time or for the inductance of the RLC analogy, $R$ and $C$ for resistance and capacitance, $p$ for a pole location, and $\tau_c$ for the time constant of the SIMC rule.

\Needspace{15\baselineskip}
\paragraph{Signals and PID (\cref{sec:realization}).}
\begin{center}\small
\begin{tabularx}{\textwidth}{@{}lX@{}}
\toprule
$r_t,y_t,e_t$ & reference, output, and error $e_t=r_t-y_t$\\
$\Delta$ & sampling period\\
$I_t$, $D_t$ & integral $I_t=I_{t-1}+\Delta e_t$ and difference quotient $(e_t-e_{t-1})/\Delta$\\
$d$ & error difference $e_t-e_{t-1}$ in the example ($d=D_t$ since $\Delta=1$)\\
$k=(K_p,K_i,K_d)^\T$, $\phi$ & PID gains and PID features $(e_t,I_t,D_t)^\T$\\
$s_t$, $A_s$, $b_s$ & two-state memory $(I_t,e_t)^\T$ and its update\\
$c=(c_1,c_2)^\T$, $c_0$ & two-state readout, $u_t=c^\T s_{t-1}+c_0e_t$\\
$\Gamma_\Delta$, $\theta$ & readout map and readout parameters $\theta=(c,c_0)=\Gamma_\Delta k$\\
$\tilde s_t$, $\tilde A$, $\tilde b$, $\tilde c$ & three-state realization \eqref{eq:threestate}\\
$a_0,a_1,a_2$ & incremental-form coefficients\\
$f,i,g$ & forget gate, input gate, and candidate of the gated variant\\
$h_t$, $c^{\rm cell}_t$, $o_t$ & LSTM hidden state of width $n_h$, cell state, and output gate (\cref{prop:lstm}); $h_t$ is also the residual's hidden state in \eqref{eq:residual}, and the bold $\mathbf o_t$ is the observation of \cref{sec:warm}\\
$W,b_o$ & output weights and bias of a bounded-output readout (\cref{prop:lstm})\\
\bottomrule
\end{tabularx}
\end{center}

\Needspace{15\baselineskip}
\paragraph{One-step objective and theorems (\cref{sec:learning,sec:improvement}).}
\begin{center}\small
\begin{tabularx}{\textwidth}{@{}lX@{}}
\toprule
$z\in\R^n$, $\nu$, $z^+$ & augmented state (plant state and controller memory), its design distribution, and the next state\\
$a(z)$, $b$, $\xi$ & drift, constant input vector, and transition noise in \eqref{eq:plant}\\
$H$, $R$, $M$ & state weight, action weight, and $M=R+b^\T Hb$\\
$\mathcal L$, $\mathcal L_{\rm irr}$ & one-step loss and its policy-independent part\\
$u_\star$ & pointwise one-step optimum $-M^{-1}b^\T Ha(z)$\\
$\Sigma$, $k_\star$ & PID feature Gram matrix and optimal PID gains\\
$\eta$, $\ell$, $q$ & step size, iteration index, and contraction factor\\
$\mathcal H_\theta$, $\kappa(\cdot)$ & Hessian in readout coordinates; condition number\\
$\psi\in\R^p$, $\beta$, $\beta_0$ & nonlinear features, their coefficients, and a constant offset\\
$C_{\psi\phi}$, $\widetilde\psi$, $G$ & regression of $\psi$ on $\phi$, residualized features, and their Gram matrix\\
$m_\phi$, $m_\psi$ & cross moments $\E[\phi u_\star]$ and $\E[\widetilde\psi u_\star]$\\
$\beta_\star$, $k_{\rm raw}$, $k_{\rm lin}$ & optimal residual coefficient; raw-coordinate PID block; linear part of $u_\star$\\
$G_\psi$ & Gram matrix of the residual features in \cref{cor:warm}: $G$, or $\E[\psi\psi^\T]$ in raw coordinates\\
$P_\phi$ & orthogonal projection onto the PID span in $L^2(\nu)$\\
$A$, $\alpha$ & example matrix and nonlinearity strength in \eqref{eq:example}\\
\bottomrule
\end{tabularx}
\end{center}

\Needspace{15\baselineskip}
\paragraph{Finite samples and long horizons (\cref{sec:long}, \cref{app:finite,app:longhorizon}).}
\begin{center}\small
\begin{tabularx}{\textwidth}{@{}lX@{}}
\toprule
$N$, $\widehat\Sigma$, $\widehat m_\phi$ & sample size and empirical moments\\
$\widehat\epsilon$ & noise term of the empirical minimizer\\
$\mu$ & lower eigenvalue bound (\cref{app:finite}) or Polyak--\L{}ojasiewicz constant (\cref{app:longhorizon})\\
$T$, $J_T$, $\vartheta$ & horizon, normalized tracking cost, and controller parameters\\
$J_T^{\PID}$ & tracking cost of the practical PID alone (\cref{prop:warm})\\
$\chi_t$, $F$, $\xi_t$, $S_t$ & closed-loop state, its update, exogenous randomness, and forward sensitivity\\
$\mathcal V$, $L$, $\mu_{\rm qg}$ & region of the iterates, smoothness constant, and quadratic-growth constant\\
\bottomrule
\end{tabularx}
\end{center}

\Needspace{15\baselineskip}
\paragraph{Practical controllers and experiments (\crefrange{sec:long}{sec:experiments}).}
\begin{center}\small
\begin{tabularx}{\textwidth}{@{}lX@{}}
\toprule
$\lambda_{\rm f}$ & derivative-filter coefficient (distinct from eigenvalues $\lambda_j$)\\
$v_t$, $k_{\rm aw}$, $\tau_{\rm t}$ & integral state in control units, anti-windup fraction, and tracking time constant\\
$u_t$, $\bar u_t$, $u_{\max}$ & commanded action, applied (clipped) action, and actuator limit\\
$\rho_t$, $\rho_{\max}$ & bounded residual and its amplitude\\
$W_h,W_x,b_h$, $W_\rho,b_\rho$ & residual network: recurrent, input and output weights and biases\\
$\mathbf o_t$ & observation of the learned controllers, $(e_t,I_t,D_t,\bar u_{t-1})$\\
$\vartheta_{\rm P},\vartheta_{\rm m},\vartheta_{\rm o}$ & PID, memory and output-layer parameters of the residual controller (\cref{sec:warm})\\
$p_t$ & sensitivity of the episode cost to an additive action perturbation at step $t$\\
$\delta_u$ & constant actuator offset in the benchmark families\\
$c_{\rm ex}$ & leading-order constant of the learned excess cost in E3\\
$s_e$, $s_u$ & error and action scales in \eqref{eq:longcost}\\
$\zeta,\omega,\gamma$ & oscillator damping, natural frequency, and Duffing coefficient\\
$a_{\rm p},b_{\rm p},c_{\rm p}$ & plant-family parameters (Hammerstein; $b_{\rm p}$ also the oscillator input gain)\\
$w$, $w_t$ & process noise of the benchmark families\\
$\mathcal R$ & relative cost against the PID envelope ($\mathcal R_f$ for family $f$)\\
$\Delta F_f$ & failure-rate difference, residual minus envelope, in family $f$\\
$c_{\rm cap}$ & per-step cost charged after an aborted rollout\\
\bottomrule
\end{tabularx}
\end{center}

\end{document}